\documentclass[10pt]{amsart}
\usepackage[T1]{fontenc}
\usepackage{geometry}
\usepackage[latin1] {inputenc}
\usepackage{amsmath}
\usepackage{amsfonts, amssymb, textcomp}
\usepackage[colorlinks=flase, linkcolor=red,urlcolor=green, citecolor=blue]{hyperref}
\usepackage{subeqnarray}

\usepackage{latexsym}
\usepackage{fancyhdr}
\usepackage{longtable}
\usepackage{amsmath, amssymb}
\usepackage{graphicx}
\usepackage{enumerate}

\numberwithin{equation}{section}

\theoremstyle{plain}
\newtheorem{proposition}{Proposition}[section]
\newtheorem{theorem}[proposition]{Theorem}
\newtheorem{lemma}[proposition]{Lemma}
\newtheorem{corollary}[proposition]{Corollary}

\newtheorem{example}[proposition]{Example}

\newtheorem{remark}[proposition]{Remark}

\newcommand{\RR}{\mathbb{R}}
\newcommand{\CC}{\mathbb{C}}
\newcommand{\NN}{\mathbb{N}}

\let\on=\operatorname

\newenvironment{demo}[1]{\par\smallskip\noindent{\bf #1.}}{\par\smallskip}

\newsavebox{\fmbox}
\newenvironment{fmpage}[1]
 {\begin{lrbox}{\fmbox}\begin{minipage}{#1}}
 {\end{minipage}\end{lrbox}\fbox{\usebox{\fmbox}}}

\title[Solid hulls and cores of weighted entire functions]
{Solid hulls and cores of classes of weighted entire functions defined in terms of associated weight functions}

\author[G.~Schindl]{Gerhard Schindl}

\address{G.~Schindl: Fakult\"at f\"ur Mathematik, Universit\"at Wien, Oskar-Morgenstern-Platz~1, A-1090 Wien, Austria.}
\email{gerhard.schindl@univie.ac.at}

\begin{document}

\begin{abstract}
In the spirit of very recent articles by J. Bonet, W. Lusky and J. Taskinen we are studying the so-called solid hulls and cores of spaces of weighted entire functions when the weights are given in terms of associated weight functions coming from weight sequences. These sequences are required to satisfy certain (standard) growth and regularity properties which are frequently arising and used in the theory of ultradifferentiable and ultraholomorphic function classes (where also the associated weight function plays a prominent role). Thanks to this additional information we are able to see which growth behavior the so-called ''Lusky-numbers'', arising in the representations of the solid hulls and cores, have to satisfy resp. if such numbers can exist.
\end{abstract}

\thanks{G. Schindl is supported by FWF-Projects J~3948-N35 and P32905}
\date{\today}

\maketitle


\section{Introduction}\label{Introduction}
Spaces of weighted entire functions are defined as follows
$$H^{\infty}_v(\CC):=\{f\in H(\CC): \|f\|_v:=\sup_{z\in\CC}|f(z)|v(|z|)<+\infty\},$$
and the weight $v:[0,+\infty)\rightarrow(0,+\infty)$ is usually assumed to be continuous, non-increasing and rapidly decreasing, i.e. $\lim_{r\rightarrow+\infty}r^kv(r)=0$ for all $k\ge 0$. In the recent publications \cite{BonetTaskinen18} and \cite{BonetLuskyTaskinen19} the authors have studied the so-called {\itshape solid hulls} and {\itshape solid cores} of such spaces, using the identification of $f(z)=\sum_{j=0}^{+\infty}a_jz^j$ with its sequence of Taylor-coefficients $(a_j)_{j\in\NN}$. For more references and historical background we refer to the introductions of these papers.\vspace{6pt}

It has turned out that the so-called regularity condition $(b)$ on the weight $v$, see \cite[$(2.2)$]{BonetTaskinen18} and \cite[Definition 2.1]{BonetLuskyTaskinen19} and \eqref{22} in the present work, plays the key-role for a more concrete description of the solid hulls and cores of $H^{\infty}_v(\CC)$. It is weaker than condition $(B)$ introduced on \cite[p. 20]{Lusky06}, see \cite[Rem. 2.7]{BonetTaskinen18}, and the arising expressions in $(b)$ have already been studied in \cite{Lusky06} as well. Verifying this condition for concrete given examples requires quite technical computations and might be challenging: The expressions under consideration are involving the extremal points of the functions $r\mapsto r^kv(r)$ and one has to find and compute a strictly increasing sequence of positive real numbers, the so-called {\itshape Lusky numbers}.\vspace{6pt}

The goal of this paper is to study the situation when $v(r)=\exp(-\omega_M(r))$, with $\omega_M$ denoting the so-called {\itshape associated weight function} (see \eqref{assofunc}), and $M=(M_p)_{p\in\NN}$ a given sequence of positive real numbers satisfying some basic regularity and growth properties.\vspace{6pt}

First, this question has been motivated by recognizing that the family of weights studied in \cite[Sect. 3]{BonetTaskinen18} corresponds (up to an equivalence of weight functions) to the associated weight functions coming from the {\itshape Gevrey sequences} $(p!^s)_{p\in\NN}$, $s>0$, arising frequently in the theory of ultradifferentiable and ultraholomorphic functions.

Second, by a known characterizing result concerning so-called ultradifferentiable operators from \cite{Komatsu73} (see Proposition \ref{Komaprop45}), classes of weighted entire functions with the weight $v(r)=\exp(-\omega_M(r))$ are naturally arising also in the ultradifferentiable setting.\vspace{6pt}

Based on these observations the main idea has then been to connect two areas of research and exploit the additional information on the (standard) growth and regularity properties of the underlying weight sequence in order to verify condition $(b)$, more precisely: Apply this approach to compute, via alternative techniques, explicitly the Lusky numbers, get knowledge about their possible growth resp. see and decide that such numbers cannot exist. This has also led us to the following questions: How are the (standard) properties on a weight sequence, which are arising frequently in the ultradifferentiable and ultraholomorphic framework, related to the regularity condition $(b)$? Do ''nice and very regular'' sequences, e.g. {\itshape strongly regular sequences} (see \cite[1.1]{Thilliezdivision}), always admit the existence of Lusky numbers?\vspace{6pt}

Summarizing, on the one hand we have been able to see several connections between the required growth properties and so the weight sequence setting helps to compute the Lusky numbers resp. to see that such numbers cannot exist. E.g. we have been able to see that too fast increasing weight sequences do not admit the existence of Lusky numbers, see Lemma \ref{alternativedeltaconsequ}. But, on the other hand, we have been able to construct (counter-)examples showing that, roughly speaking, the required resp. desired notions of regularity in the ultradifferentiable  and ultraholomorphic world and in the weighted entire world fall apart.\vspace{6pt}

We summarize now the content of this article.\vspace{6pt}

After recalling and collecting all necessary basic definitions of weight sequences and their growth properties in Section \ref{sequencesdef}, we are rewriting a known characterizing result of ultradifferentiable operators in terms of solid hulls and cores, see Section \ref{Udos}.

In Section \ref{solidhullsandcores} we provide a deep study of the regularity condition $(b)$ in the weight sequence setting, see Lemma \ref{ABexpressions}, which enables us to reformulate the main results of the solid hull from \cite{BonetTaskinen18} and the solid core from \cite{BonetLuskyTaskinen19} involving the ''Lusky numbers'', see Theorems \ref{Thm25} and \ref{Thm24}.\vspace{6pt}

Then, in Section \ref{rinterpolatingsect} we study the behavior of the Lusky numbers moving from $M$ to its so-called $r$-interpolating sequence $P^{M,r}$, which has been used to prove extension results in the ultraholomorphic setting, see \cite{Schmetsvaldivia00}, \cite{injsurj}. It turns out that this natural construction can be used to determine in Theorem \ref{Pramification} the solid hulls and cores of spaces defined in terms of ramified weights, see \eqref{ramifiedweight}.

Using the derived formulas for the condition $(b)$, in Section \ref{firstcounter} we are able to give examples of weight sequences $M$ which have many good growth and regularity properties but do not admit the existence of Lusky numbers.

A first new example, related to the weight $v(r)=\exp(-(\log(1+r))^2)$ which is given by the so-called $q$-Gevrey sequences $(q^{p^2})_{p\in\NN}$, $q>1$, is studied in detail in Section \ref{qGevrey}. We are able to compute the Lusky numbers, obtain closed explicit expressions for the solid hull and core and can prove that the space is solid, see Corollaries \ref{Thm25qgevrey} and \ref{Thm25qgevrey1}.\vspace{6pt}

In Section \ref{Alternativesequencedelta}, by using an auxiliary sequence $(\delta_p)_p$, which is measuring the growth of the quotients $\mu_p:=M_p/M_{p-1}$, we are able to provide a more detailed and precise study of the connection between the growth of $M$ and the (possibly) existing Lusky numbers.

On the one hand, in Corollary \ref{alternativedeltaconsequ1cor} we derive necessary growth restrictions for sequences of integers which can serve as Lusky numbers in the weight sequence setting, see \eqref{Luskynecessary}. Conversely, in Proposition \ref{ajexample} we show that for each sequence of integers $(a_j)_j$ satisfying these growth restrictions, we can associate a weight sequence $M$ such that $(a_j)_j$ can be used as Lusky numbers for this particular constructed sequence. So, roughly speaking, the mapping $M\mapsto(a_j)_j$ is surjective and the weight sequence setting is in this sense sufficiently large enough.\vspace{6pt}

In Lemma \ref{whatcanhappen1} we construct a (counter-)example of a sequence being equivalent to a strongly regular sequence but not admitting Lusky numbers which underlines the different behavior of condition $(b)$ compared with the notion of regularity in the ultradifferentiable and ultraholomorphic setting.\vspace{6pt}

In Section \ref{Gevreyexample}, as a second concrete example, we compute the Lusky numbers for the Gevrey sequences and are proving the characterization given in \cite[Theorem 3.1]{BonetTaskinen18} by completely different methods.\vspace{6pt}

Section \ref{fromfcttosequ} provides some information about the idea when starting with a given abstract weight function $v$, satisfying the technical condition \eqref{vconvexity}, and then considering the associated weight sequence $M^v$ analogously as it has been done in the ultradifferentiable setting. The motivation behind this approach is that it shows how the results obtained in this article can help to get information for abstractly given weight functions $v$ as well.\vspace{6pt}

The solid hulls and solid cores of weighted holomorphic functions on the unit disk have been studied in \cite{BonetTaskinendisc18} and \cite{BonetLuskyTaskinen19} and it has turned out that also in this situation the regularity condition $(b)$ becomes crucial for a precise description, again in terms of the Lusky numbers. For the sake of completeness, in Section \ref{disksection} we investigate this notion also in the weight sequence setting, see Theorem \ref{Thm2524disc}. Unfortunately, here the arising expressions are becoming much more involved (see Lemma \ref{ABexpressionsdisc}) and concrete computations for the Lusky numbers seem to be much more complicated and involved.\vspace{6pt}

\textbf{Acknowledgements.} The author of this article thanks the anonymous referees for their careful reading and their valuable suggestions which have improved and clarified the presentation of the results. The author wishes to thank Prof. Jos\'{e} Bonet Solves from the Universitat Polit\`{e}cnica de Val\`{e}ncia and Prof. Wolfgang Lusky from the University of Paderborn for fruitful and interesting discussions and for giving helping and clarifying explanations during the preparation of this work. Moreover, he wants to thank Armin Rainer from the University of Vienna for giving some suggestions during reading a preliminary version. Finally, the author expresses his deepest thank to his friend David N. Nenning, also from the University of Vienna, for carefully reading a preliminary version of this article and for giving several important hints and ideas.







\section{Notation}\label{sequencesdef}
\subsection{General notation}
We write $\NN:=\{0,1,2,\dots\}$ and $\NN_{>0}:=\{1,2,3,\dots\}$. With $\lfloor x\rfloor$ we denote the integer part of any given $x>0$, with $\lceil x\rceil$ the smallest integer greater or equal than $x$. $H(\CC)$ denotes the space of entire functions.

\subsection{Weight sequences}
Given a sequence $M=(M_j)_j\in\RR_{>0}^{\NN}$ we also use $\mu_j:=\frac{M_j}{M_{j-1}}$, $\mu_0:=1$. $M$ is called {\itshape normalized} if $1=M_0\le M_1$ holds true which can always be assumed without loss of generality. For any $s>0$ we write $M^s=(M^s_j)_{j\in\NN}$ for the $s$-th power of $M$.\vspace{6pt}

In the following we collect several growth and regularity properties for $M$ which will be used later on. These conditions arise frequently and are standard in the ultradifferentiable and ultraholomorphic weight sequence setting (e.g. see \cite{Komatsu73}).

$M$ is called {\itshape log-convex} if
$$\forall\;j\in\NN_{>0}:\;M_j^2\le M_{j-1} M_{j+1},$$
equivalently if $(\mu_j)_{j\ge 1}$ is non-decreasing. If $M$ is log-convex and normalized, then both $j\mapsto M_j$ and $j\mapsto(M_j)^{1/j}$ are non-decreasing and $(M_j)^{1/j}\le\mu_j$ for all $j\in\NN_{>0}$ (e.g. see \cite[Lemma 2.0.4]{diploma}). If the sequence $m:=\left(\frac{M_j}{j!}\right)_j$ is log-convex, then $M$ is called {\itshape strongly log-convex}, denoted by \hypertarget{slc}{$(\text{slc})$}.

For our purpose it is convenient to consider sequences belonging to the set
$$\hypertarget{LCset}{\mathcal{LC}}:=\{M\in\RR_{>0}^{\NN}:\;M\;\text{is normalized, log-convex},\;\lim_{j\rightarrow+\infty}(M_j)^{1/j}=+\infty\}.$$
We see that $M\in\hyperlink{LCset}{\mathcal{LC}}$ if and only if $1=\mu_0\le\mu_1\le\dots$, $\lim_{j\rightarrow+\infty}\mu_j=+\infty$ (e.g. see \cite[p. 104]{compositionpaper}) and there is a one-to-one correspondence between $M$ and $\mu=(\mu_j)_j$ by taking $M_p:=\prod_{i=0}^p\mu_i$.\vspace{6pt}

In the ultradifferentiable and ultraholomorphic setting the following conditions on $M$ arise frequently in the literature.

We say that $M$ has {\itshape derivation closedness}, denoted by \hypertarget{dc}{$(\text{dc})$}, if
$$\exists\;D\ge 1\;\forall\;j\in\NN:\;M_{j+1}\le D^{j+1} M_j\Longleftrightarrow\mu_{j+1}\le D^{j+1},$$
and $M$ has the stronger condition {\itshape moderate growth}, denoted by \hypertarget{mg}{$(\text{mg})$}, if
$$\exists\;C\ge 1\;\forall\;j,k\in\NN:\;M_{j+k}\le C^{j+k} M_j M_k.$$
It is known (e.g. see \cite[Lemma 2.2]{whitneyextensionweightmatrix}) that for any given $M\in\hyperlink{LCset}{\mathcal{LC}}$ condition \hyperlink{mg}{$(\on{mg})$} is equivalent to having $\sup_{j\in\NN}\frac{\mu_{2j}}{\mu_j}<+\infty$.

We say that $M$ has \hypertarget{beta3}{$(\beta_3)$} if
$$\exists\;Q\in\NN_{\ge 2}:\;\;\;\liminf_{j\rightarrow\infty}\frac{\mu_{Qj}}{\mu_j}>1,$$
(see \cite{BonetMeiseMelikhov07}) and $M$ has the stronger condition \hypertarget{beta1}{$(\beta_1)$} (introduced in \cite{petzsche}), if
$$\exists\;Q\in\NN_{\ge 2}:\;\;\;\liminf_{j\rightarrow+\infty}\frac{\mu_{Qj}}{\mu_j}>Q.$$
$M$ has \hypertarget{gammar}{$(\gamma_1)$} if
$$\sup_{j\in\NN_{>0}}\frac{\mu_j}{j}\sum_{k\ge j}\frac{1}{\mu_j}<+\infty.$$
In the literature \hyperlink{gammar}{$(\gamma_1)$} is also called ``strong nonquasianalyticity condition'', see \cite{petzsche} and  \cite{Komatsu73} and for any $M\in\hyperlink{LCset}{\mathcal{LC}}$ in \cite{petzsche} it has been shown that $\hyperlink{beta1}{(\beta_1)}\Longleftrightarrow\hyperlink{gammar}{(\gamma_1)}$.

A smaller class than \hyperlink{LCset}{$\mathcal{LC}$} is the set $\mathcal{SR}$ defined by\vspace{6pt}

\centerline{$M\in\hypertarget{SRset}{\mathcal{SR}}$, if $M$ is normalized and has \hyperlink{slc}{$(\text{slc})$}, \hyperlink{mg}{$(\text{mg})$} and \hyperlink{gamma1}{$(\gamma_1)$}.}\vspace{6pt}

Using this notation we see that $M\in\hyperlink{SRset}{\mathcal{SR}}$ if and only if $m$ is a {\itshape strongly regular sequence} in the sense of \cite[1.1]{Thilliezdivision}.\vspace{6pt}

Let $M,N\in\RR_{>0}^{\NN}$ be given, we write $M\hypertarget{preceq}{\preceq}N$ if $\sup_{j\in\NN_{>0}}\left(\frac{M_j}{N_j}\right)^{1/j}<+\infty$. We call $M$ and $N$ {\itshape equivalent}, denoted by $M\hypertarget{approx}{\approx}N$, if $M\hyperlink{preceq}{\preceq}N$ and $N\hyperlink{preceq}{\preceq}M$. This relation is characterizing for $M,N\in\hyperlink{LCset}{\mathcal{LC}}$ the equivalence of ultradifferentiable function classes (e.g. see \cite[Prop. 2.12]{compositionpaper}).

We mention that in \cite[Prop. 1.1]{petzsche} it has been shown that \hyperlink{gammar}{$(\gamma_1)$} for log-convex $M$ implies that there does exist an equivalent sequence $N$ having \hyperlink{slc}{$(\text{slc})$}, so \hyperlink{gammar}{$(\gamma_1)$} ''implies'' \hyperlink{slc}{$(\text{slc})$}.\vspace{6pt}

A prominent example is $G^s:=(j!^s)_{j\in\NN}$, the so-called {\itshape Gevrey-sequence} of index $s>0$. If $s>1$, then it is straightforward to check that $G^s\in\hyperlink{SRset}{\mathcal{SR}}$ (e.g. see again \cite[1.1]{Thilliezdivision}).\vspace{6pt}

Let $M\in\RR_{>0}^{\NN}$ (with $M_0=1$) be given. Then the
associated function $\omega_M: \RR_{\ge 0}\rightarrow\RR\cup\{+\infty\}$ is defined by
\begin{equation}\label{assofunc}
\omega_M(t):=\sup_{j\in\NN}\log\left(\frac{t^j}{M_j}\right)\,\,\,\text{
	for }\,t>0,\quad \omega_M(0):=0.
\end{equation}

If $\liminf_{j\rightarrow+\infty}(M_j)^{1/j}>0$, then $\omega_M(t)=0$ for sufficiently small $t$, since $\log\left(\frac{t^j}{M_j}\right)<0\Leftrightarrow t<(M_j)^{1/j}$ holds for all $j\in\NN_{>0}$. Moreover under this assumption $t\mapsto\omega_M(t)$ is a continuous nondecreasing function, which is convex in the variable $\log(t)$ and tends faster to infinity than any $\log(t^j)$, $j\ge 1$, as $t\rightarrow+\infty$. $\lim_{j\rightarrow+\infty}(M_j)^{1/j}=+\infty$ implies that $\omega_M(t)<+\infty$ for any finite $t$ which shall be considered as a basic assumption for defining $\omega_M$. We refer to \cite[Chapitre I]{mandelbrojtbook} and \cite[Def. 3.1]{Komatsu73}.\vspace{6pt}

Given $M\in\hyperlink{LCset}{\mathcal{LC}}$, then by \cite[1.8 III]{mandelbrojtbook} we get: One has $\omega_M(t)=0$ on $[0,\mu_1]$ and for all $t\ge\mu_1(\ge 1)$, when $j_t\in\NN_{>0}$ is denoting the index such that $\mu_{j_t}\le t<\mu_{j_t+1}$ is valid, then we get
\begin{equation}\label{assovsmu}
\omega_M(t)=j_t\log(t)-\log(M_{j_t}).
\end{equation}
Note that for $M\in\hyperlink{LCset}{\mathcal{LC}}$ we have $\lim_{j\rightarrow+\infty}\mu_j=+\infty$ (e.g. see \cite[p. 104]{compositionpaper}).\vspace{6pt}

If $M,N\in\RR_{>0}^{\NN}$ (with $M_0=N_0=1$), then $M\hyperlink{preceq}{\preceq}N$ does imply $\omega_N(t)\le\omega_M(Ct)$ for some $C\ge 1$ and all $t\ge 0$. If $M,N\in\hyperlink{LCset}{\mathcal{LC}}$, then the converse holds true as well, see \cite[Lemma 3.8]{Komatsu73}.

Given two (associated weight) functions $\sigma,\tau$, we write $\sigma\hypertarget{ompreceq}{\preceq}\tau$ if
$$\tau(t)=O(\sigma(t))\;\text{as}\;t\rightarrow+\infty$$
and call them equivalent, denoted by $\sigma\hypertarget{sim}{\sim}\tau$, if
$$\sigma\hyperlink{ompreceq}{\preceq}\tau\;\text{and}\;\tau\hyperlink{ompreceq}{\preceq}\sigma.$$

It is known that the equivalence of ultradifferentiable function classes defined by (Braun-Meise-Taylor) weight functions is characterized by this relation (see \cite[Cor. 5.17]{compositionpaper}) and for any $s>0$ the mapping $t\mapsto t^s$ is equivalent to $\omega_{G^{1/s}}$.

\section{Solid hulls, cores and ultradifferentiable operators}\label{Udos}
We recall now briefly the notion of {\itshape solid sub-
	and superspaces} for spaces of (complex) sequences, e.g. see
\cite{andersonshields}. Let $A$ be a vector space of sequences, then
$A$ is said to be {\itshape solid} if $(a_j)_j\in A$ does imply
$(b_j)_j\in A$ for all sequences satisfying $|b_j|\le|a_j|$,
$\forall j\in\NN$.

In \cite[Lemma 2]{andersonshields} it has been shown that for any
given sequence space $A$ there does exist $s(A)$, the {\itshape
	largest solid subspace (or solid core)} of $A$, and there does exist
$S(A)$, the {\itshape smallest solid superspace (or solid hull)}, of
$A$. We have
\begin{equation*}\label{solidcore}
s(A)=\{(b_j)_{j\in \NN}: (b_j \lambda_j)_{j\in \NN}\in A,
\forall\,(\lambda_j)_{j\in \NN}\in l^{\infty}\}
\end{equation*} and
\begin{equation*}\label{solidalternative}
S(A)=\{(b_j)_{j\in \NN}:\,\exists\,(a_j)_{j\in \NN}\in A:\,|b_j|\le|a_j|,\,\forall\,j\in\NN\},
\end{equation*}
e.g. see \cite[p. 594]{BonetTaskinen18}. It is clear that $A\subseteq B$ does imply $s(A)\subseteq s(B)$ and $S(A)\subseteq S(B)$.\vspace{6pt}

We start by recalling \cite[Prop. 4.5]{Komatsu73} where the following characterization has been shown and this has been a main motivation for writing this article:

\begin{proposition}\label{Komaprop45}
Let $M\in\hyperlink{LCset}{\mathcal{LC}}$ and an entire function $P(\xi)=\sum_{j=0}^{+\infty}a_j\xi^j$ be given.
\begin{itemize}
\item[$(i)$] The following are equivalent (''Roumieu-type variant''):
\begin{itemize}
\item[$(a)$] $\exists\;L,C>0\;\forall\;j\in\NN:\:\;\;|a_j|\le\frac{CL^j}{M_j}$,
\item[$(b)$] $\exists\;L,C>0\;\forall\;\xi\in\CC:\:\;\;|P(\xi)|\le C\exp(\omega_M(L\xi))$.
\end{itemize}
\item[$(ii)$]
The following are equivalent (''Beurling-type variant''):
\begin{itemize}
\item[$(a)$] $\forall\;L>0\;\exists\;C>0\;\forall\;j\in\NN:\:\;\;|a_j|\le\frac{CL^j}{M_j}$,
\item[$(b)$] $\forall\;L>0\;\exists\;C>0\;\forall\;\xi\in\CC:\:\;\;|P(\xi)|\le C\exp(\omega_M(L\xi))$.
\end{itemize}
\end{itemize}
In both cases $\omega_M(L\xi)$ means $\omega_M(L|\xi|)$ (radial extension to $\CC$). The proof shows that in $(a)\Rightarrow(b)$ we have to replace $L$ by $2L$, in $(b)\Rightarrow(a)$ we can take the same constant $L$ (since we are considering the case dimension $d=1$).
\end{proposition}

An operator of the form $P(\partial):=\sum_{j=0}^{+\infty}a_j\partial^j$ with $(a_j)_j$ satisfying $(i)(a)$ resp. $(ii)(a)$ above is called an {\itshape ultradifferential operator of Roumieu- resp. Beurling-type.} Hence the previous result motivates also from the point of view of studying problems in the ultradifferentiable setting to consider spaces of weighted entire functions defined as follows:

For any $M\in\hyperlink{LCset}{\mathcal{LC}}$ and $c>0$ we set
\begin{equation}\label{weights}
v_{M,c}(r):=\exp(-\omega_M(cr)),\;\;\;r\in[0,+\infty),
\end{equation}
and
\begin{align*}
H^{\infty}_{v_{M,c}}(\CC)&:=\{f\in H(\CC): \|f\|_{v_{M,c}}:=\sup_{z\in\CC}|f(z)|v_{M,c}(|z|)<+\infty\}
\\&
=\{f\in H(\CC): \|f\|_{v_{M,c}}:=\sup_{z\in\CC}|f(z)|\exp(-\omega_M(c|z|))<+\infty\}.
\end{align*}
For $c=1$ we will write $v_M$ instead of $v_{M,1}$. We call $v:[0,+\infty)\rightarrow(0,+\infty)$ a {\itshape weight function}, if $v$ is {\itshape continuous, non-increasing and rapidly decreasing,} i.e. $\lim_{r\rightarrow+\infty}r^kv(r)=0$ for all $k\ge 0$. This is the same notion to be a weight function as it has been considered in \cite{BonetTaskinen18} and \cite{BonetLuskyTaskinen19}.

First note that by the properties of $\omega_M$, by definition each $v_{M,c}$ is a weight function in the previous sense (with having $v_{M,c}(0)=1$), see again \cite[Chapitre I]{mandelbrojtbook} and \cite[Def. 3.1]{Komatsu73}.\vspace{6pt}

Given $M,N\in\hyperlink{LCset}{\mathcal{LC}}$ with $M\hyperlink{approx}{\approx}N$, then for some $D\ge 1$ and all $t\ge0$ we have $\omega_N(D^{-1}t)\le\omega_M(t)\le\omega_N(Dt)$, which implies
\begin{equation}\label{Hinfinclusion}
\forall\;c>0:\;\;\;H^{\infty}_{v_{N,c}}(\CC)\subseteq H^{\infty}_{v_{M,Dc}}(\CC)\subseteq H^{\infty}_{v_{N,D^2c}}(\CC).
\end{equation}
If we set for $c>0$
$$U_{M,c}:=\left\{(a_j)_{j\in\NN}\in\CC^{\NN}:\;\exists\;D>0\;\forall\;j\in\NN:\:\;\;|a_j|\le\frac{Dc^j}{M_j}\right\},$$
then Proposition \ref{Komaprop45} tells us that the sequence of Taylor coefficients $(a_j)_{j\in\NN}\in\CC^{\NN}$ of any $f(z)=\sum_{j=0}^{+\infty}a_jz^j\in H^{\infty}_{v_{M,c}}(\CC)$ has to belong to $U_{M,c}$ and any sequence $(a_j)_{j\in\NN}\in U_{M,c}$ yields an entire $f(z)=\sum_{j=0}^{+\infty}a_jz^j\in H^{\infty}_{v_{M,2c}}(\CC)$.\vspace{6pt}

By definition it is clear that $U_{M,c}$ is solid, i.e. $U_{M,c}=S(U_{M,c})=s(U_{M,c})$, hence Proposition \ref{Komaprop45} tells us that (by identifying a function $f(z)=\sum_{j=0}^{+\infty}a_jz^j\in H^{\infty}_{v_{M,c}}(\CC)$ with its sequence of Taylor coefficients $(a_j)_j$):
$$\forall\;c>0:\;\;\;U_{M,c/2}\subseteq s(H^{\infty}_{v_{M,c}}(\CC))\subseteq H^{\infty}_{v_{M,c}}(\CC)\subseteq S(H^{\infty}_{v_{M,c}}(\CC))\subseteq U_{M,c}.$$

In \cite[Proposition 1.1]{BonetTaskinen18} a first characterization of the solid core $s(H^{\infty}_{v}(\CC))$ has been obtained and this results takes in our setting the following form:

\begin{proposition}\label{Prop11}
Let $M\in\hyperlink{LCset}{\mathcal{LC}}$ be given and $c>0$. Then the solid core of $H^{\infty}_{v_{M,c}}$ is given by
$$s(H^{\infty}_{v_{M,c}}(\CC))=\{(b_j)_{j\in\NN}\in\CC^{\NN}: \sup_{r>0}v_{M,c}(r)\sum_{j=0}^{+\infty}|b_j|r^j<+\infty\}.$$
\end{proposition}

\section{Solid hulls and cores of weighted entire functions}\label{solidhullsandcores}
In order to get more precise information on the solid core $s(H^{\infty}_{v_{M,c}}(\CC))$ and the solid hull $S(H^{\infty}_{v_{M,c}}(\CC))$ we have to study the regularity condition $(b)$, see \cite[$(2.2)$]{BonetTaskinen18} and \cite[Definition 2.1]{BonetLuskyTaskinen19}.

Recall that, as commented above, the weight $z\mapsto\exp(-|z|^s)$, $s>0$, is corresponding (up to an equivalence) to $v_{G^{1/s},1}$, with $G^{1/s}$ denoting the {\itshape Gevrey sequence} with index $1/s$.

\subsection{On the regularity condition $(b)$ in the weight sequence setting}\label{regularitybcond}

For any $k\ge 0$, $c>0$ and $r\ge 0$ we set $G^k_{M,c}(r):=r^kv_{M,c}(r)$, more precisely one has:
\begin{align*}
G^k_{M,c}(r)&=r^k\exp(-\omega_M(cr))=r^k\exp(-\sup_{p\in\NN}\log((cr)^p/M_p))=r^k\exp(\inf_{p\in\NN}-\log((cr)^p/M_p))
\\&
=r^k\exp(\inf_{p\in\NN}\log(M_p/(cr)^p))=r^k\inf_{p\in\NN}\frac{M_p}{(cr)^p}=r^kh_M\left(\frac{1}{cr}\right),
\end{align*}
with $h_M(t):=\inf_{p\in\NN}t^pM_p$, $t>0$, $h_M(0):=0$. $h_M$ is denoting another auxiliary function arising in extension theorems in the ultradifferentiable setting, e.g. see \cite{Dynkin80} or \cite{ChaumatChollet94}. The connection to $\omega_M$ is given by $h_M(t)=\exp(-\omega_M(1/t))$, hence $h_M$ is continuous and non-decreasing and $h_M(t)=1$ for all $t>0$ sufficiently large.

Since $M\in\hyperlink{LCset}{\mathcal{LC}}$, by \eqref{assovsmu} we have $h_M(\frac{1}{cr})=1$ if $\frac{1}{cr}\ge\frac{1}{\mu_1}\Leftrightarrow r\le\frac{\mu_1}{c}$ and for all $p\ge 1$: $$h_M\left(\frac{1}{cr}\right)=\frac{M_p}{(cr)^p}\;\;\;\text{if}\;\frac{1}{\mu_{p+1}}\le\frac{1}{cr}<\frac{1}{\mu_{p}}\Longleftrightarrow\mu_p< cr\le\mu_{p+1}\Longleftrightarrow\frac{\mu_p}{c}<r\le\frac{\mu_{p+1}}{c}.$$
Since $h_M$ is continuous, also $G^k_{M,c}$ is so. With $r_{k,c}$ we are denoting the global maximum point of the function $r\mapsto G^k_{M,c}(r)$ on $[0,+\infty)$.

For all values $r$ with $\mu_p< cr\le\mu_{p+1}\Leftrightarrow\frac{\mu_p}{c}<r\le\frac{\mu_{p+1}}{c}$, $p\in\NN_{\ge 1}$, we have $r^kh_M(\frac{1}{cr})=r^{k-p}\frac{M_p}{c^p}$ and so on such an interval we get $$G^k_{M,c}(r)=r^kh_M\left(\frac{1}{cr}\right)=r^{k-p}\frac{M_p}{c^p}.$$\vspace{6pt}

Next we introduce the following expressions for arbitrary $0\le k<l$, analogously to \cite[$(2.1)$]{BonetTaskinen18}:

\begin{equation}\label{ABexpressions0}
A_{M,c}(k,l):=\left(\frac{r_{k,c}}{r_{l,c}}\right)^k\frac{v_{M,c}(r_{k,c})}{v_{M,c}(r_{l,c})},\hspace{30pt}B_{M,c}(k,l):=\left(\frac{r_{l,c}}{r_{k,c}}\right)^l\frac{v_{M,c}(r_{l,c})}{v_{M,c}(r_{k,c})}.
\end{equation}

\begin{lemma}\label{ABexpressions}
Let $M\in\hyperlink{LCset}{\mathcal{LC}}$ be given and $0\le k<l$ (real numbers). Then for any $c>0$ we have that
$$A_{M,c}(k,l)=\left(\frac{\mu_{\lfloor k\rfloor+1}}{\mu_{\lfloor l\rfloor+1}}\right)^k\frac{M_{\lfloor k\rfloor+1}(\mu_{\lfloor l\rfloor+1})^{\lfloor l\rfloor+1}}{M_{\lfloor l\rfloor+1}(\mu_{\lfloor k\rfloor+1})^{\lfloor k\rfloor+1}},\hspace{20pt}B_{M,c}(k,l)=\left(\frac{\mu_{\lfloor l\rfloor+1}}{\mu_{\lfloor k\rfloor+1}}\right)^l\frac{M_{\lfloor l\rfloor+1}(\mu_{\lfloor k\rfloor+1})^{\lfloor k\rfloor+1}}{M_{\lfloor k\rfloor+1}(\mu_{\lfloor l\rfloor+1})^{\lfloor l\rfloor+1}}.$$
\end{lemma}

\demo{Proof}
Let $0\le k<l$ and $c>0$ be arbitrary, but from now on fixed.

First, if $k=0$, then $G^0_{M,c}(r)=h_M(\frac{1}{cr})$ and so the global maximum of $G^0_{M,c}$ is attained at any $r>0$ with $0<r\le\frac{\mu_1}{c}$. In this case the maximum value of $G^0_{M,c}$ is equal to $1$: Because $\omega_M$ is vanishing on $[0,\mu_1]$ we have $h_M(\frac{1}{cr})=\exp(-\omega_M(cr))=1$ for $r$ satisfying $0<r\le\frac{\mu_1}{c}$. Similarly, if $0<k<1$, then $r\mapsto r^{k-p}\frac{M_p}{c^p}$ is strictly decreasing for all $p\in\NN_{\ge 1}$. Hence by continuity of $G^k_{M,c}$, the maximum of $r\mapsto G^k_{M,c}(r)$ is attained at the right end point of $[0,\frac{\mu_1}{c}]$, because $h_M\equiv 1$ on this interval and $r\mapsto r^k$ is strictly increasing there. Summarizing, for all $0\le k<1$ we have shown $r_{k,c}=\frac{\mu_1}{c}$.\vspace{6pt}

Second, if $k\ge 1$, then by $G^k_{M,c}(r)=r^kh_M(\frac{1}{cr})=r^{k-p}\frac{M_p}{c^p}$ for all $r$ satisfying $\frac{\mu_p}{c}\le r<\frac{\mu_{p+1}}{c}$, $p\in\NN_{\ge 1}$, on such an interval the map $r\mapsto r^{k-p}\frac{M_p}{c^p}$ is strictly increasing for all $p\in\NN$ with $p<k$ and strictly decreasing for all $p>k$ (and also strictly increasing on $[0,\frac{\mu_1}{c}]$ as explained before). In the case $k\in\NN$ it is constant $\frac{M_p}{c^p}$ if $k=p$. This means that $r_{k,c}=\frac{\mu_{\lfloor k\rfloor+1}}{c}$ and if $k\in\NN$, then each $r$ with $\frac{\mu_k}{c}\le r\le\frac{\mu_{k+1}}{c}$ can be considered as a maximum value point of $G^k_{M,c}$ (with maximum value $\frac{M_k}{c^k}$). If for given $k\in\NN_{\ge 1}$ we have $\mu_k=\mu_{k+1}$, then the maximum value point coincides with $\frac{\mu_k}{c}$. (However, one can prove that w.l.o.g. we can always assume that $p\mapsto\mu_p$ is strictly increasing by changing to an equivalent sequence.)\vspace{6pt}

Thus, in the notation of \cite{BonetTaskinen18}, we can write
$$r_{k,c}=\frac{\mu_{\lfloor k\rfloor+1}}{c},\;\;\;k\ge 1,\hspace{30pt}r_{k,c}=\frac{\mu_1}{c}=\frac{\mu_{\lfloor k\rfloor+1}}{c},\;\;\;0\le k<1,$$
which should be compared for $\mu_k=k^s$, $s>1$, with \cite[p. 596]{BonetTaskinen18} (the quotient $k^s$ corresponds to the Gevrey sequence $G^{s}$).\vspace{6pt}

Since $v_{M,c}(\mu_{\lfloor k\rfloor+1}/c)=\exp(-\omega_M(\mu_{\lfloor k\rfloor+1}))=\exp(-\log(\mu_{\lfloor k\rfloor+1}^{\lfloor k\rfloor+1}/M_{\lfloor k\rfloor+1}))=\frac{M_{\lfloor k\rfloor+1}}{(\mu_{\lfloor k\rfloor+1})^{\lfloor k\rfloor+1}}$ (see \cite[1.8. III]{mandelbrojtbook}) we have
$$A_{M,c}(k,l)=\left(\frac{r_{k,c}}{r_{l,c}}\right)^k\frac{v_{M,c}(r_{k,c})}{v_{M,c}(r_{l,c})}=\left(\frac{\mu_{\lfloor k\rfloor+1}}{\mu_{\lfloor l\rfloor+1}}\right)^k\frac{v_{M,c}(\mu_{\lfloor k\rfloor+1}/c)}{v_{M,c}(\mu_{\lfloor l\rfloor+1}/c)}=\left(\frac{\mu_{\lfloor k\rfloor+1}}{\mu_{\lfloor l\rfloor+1}}\right)^k\frac{M_{\lfloor k\rfloor+1}(\mu_{\lfloor l\rfloor+1})^{\lfloor l\rfloor+1}}{M_{\lfloor l\rfloor+1}(\mu_{\lfloor k\rfloor+1})^{\lfloor k\rfloor+1}},$$
and similarly
$$B_{M,c}(k,l)=\left(\frac{r_{l,c}}{r_{k,c}}\right)^l\frac{v_{M,c}(r_{l,c})}{v_{M,c}(r_{k,c})}=\left(\frac{\mu_{\lfloor l\rfloor+1}}{\mu_{\lfloor k\rfloor+1}}\right)^l\frac{v_{M,c}(\mu_{\lfloor l\rfloor+1}/c)}{v_{M,c}(\mu_{\lfloor k\rfloor+1}/c)}=\left(\frac{\mu_{\lfloor l\rfloor+1}}{\mu_{\lfloor k\rfloor+1}}\right)^l\frac{M_{\lfloor l\rfloor+1}(\mu_{\lfloor k\rfloor+1})^{\lfloor k\rfloor+1}}{M_{\lfloor k\rfloor+1}(\mu_{\lfloor l\rfloor+1})^{\lfloor l\rfloor+1}}.$$
\qed\enddemo

In the weight sequence setting it will be convenient to assume that the numbers $k$ and $l$ are integers, say $k=a_j<a_{j+1}=l$, and then we can simplify the arising expressions to
\begin{equation}\label{Ainteresting}
A_M(a_j,a_{j+1})=A_{M,c}(a_j,a_{j+1})=\frac{(\mu_{a_{j+1}+1})^{a_{j+1}-a_j}}{\mu_{a_j+1}\cdots\mu_{a_{j+1}}}
\end{equation}
and
\begin{equation}\label{Binteresting}
B_M(a_j,a_{j+1})=B_{M,c}(a_j,a_{j+1})=\frac{\mu_{a_j+2}\cdots\mu_{a_{j+1}}}{(\mu_{a_j+1})^{a_{j+1}-a_j-1}}=\frac{\mu_{a_j+1}\mu_{a_j+2}\cdots\mu_{a_{j+1}}}{(\mu_{a_j+1})^{a_{j+1}-a_j}}.
\end{equation}

One shall note that both expressions are not depending on the given parameter $c>0$ anymore, which justifies the notation $A_M(a_j,a_{j+1})$ and $B_M(a_j,a_{j+1})$. For this recall that for any given $M\in\hyperlink{LCset}{\mathcal{LC}}$ and $c>0$ we have $\omega_N(t)=\omega_M(ct)$ with $N_k:=\frac{M_k}{c^k}$. Hence $N\hyperlink{approx}{\approx}M$ follows and more precisely $\nu_k=\frac{\mu_k}{c}$ which immediately implies that the parameter $c>0$ is cancelling out.

\begin{remark}\label{ABexpressionsremark}
Due to the discrete behavior of the weight sequence setting, we have some freedom when studying the expressions $A_M(k,l)$ and $B_M(k,l)$. As commented in the proof of Lemma \ref{ABexpressions}, for integers $k\in\NN_{>0}$ with $r_{k,c}:=\frac{r_k}{c}$ any choice $r_k\in[\mu_k,\mu_{k+1}]$ can be used as a maximal point. Then, again by \cite[1.8. III]{mandelbrojtbook}, we get for any integers $1\le k<l$:
$$\left(\frac{r_{k,c}}{r_{l,c}}\right)^k\frac{v_{M,c}(r_{k,c})}{v_{M,c}(r_{l,c})}=\left(\frac{r_k}{r_l}\right)^k\frac{M_k}{r_k^k}\frac{r_l^l}{M_l}=\frac{r_l^{l-k}}{\mu_{k+1}\cdots\mu_l},\hspace{5pt}\left(\frac{r_{l,c}}{r_{k,c}}\right)^l\frac{v_{M,c}(r_{l,c})}{v_{M,c}(r_{k,c})}=\left(\frac{r_l}{r_k}\right)^l\frac{r_k^k}{M_k}\frac{M_l}{r_l^l}=\frac{\mu_{k+1}\cdots\mu_l}{r_k^{l-k}},$$
and again there is no dependence on $c>0$ anymore.
\end{remark}
\vspace{6pt}

With this preparation, in our setting the main result \cite[Theorem 2.5]{BonetTaskinen18} takes the following form:

\begin{theorem}\label{Thm25}
	Let $M\in\hyperlink{LCset}{\mathcal{LC}}$ be given. Assume that there exists a strictly increasing sequence (of integers) $(a_j)_{j\in\NN_{\ge 1}}$, also called the ''Lusky numbers'', and constants $b$ and $K$ with $K\ge b>2$ such that
\begin{equation}\label{Thm25equ}
	b\le\min\{A_M(a_j,a_{j+1}),B_M(a_j,a_{j+1})\}\le\max\{A_M(a_j,a_{j+1}),B_M(a_j,a_{j+1})\}\le K,
\end{equation}
i.e. the regularity condition $(b)$ holds true. Then the solid hull of $H^{\infty}_{v_{M,c}}(\CC)$ is given by
\begin{equation}\label{Thm25equ1}
	S(H^{\infty}_{v_{M,c}}(\CC))=\left\{(b_j)_{j\in\NN}\in\CC^{\NN}: \sup_{j\in\NN_{\ge 1}}v_{M,c}(r_{a_j,c})\left(\sum_{a_j<l\le a_{j+1}}|b_l|^2(r_{a_j,c})^{2l}\right)^{1/2}<+\infty\right\},
\end{equation}
	or equivalently by
\begin{equation}\label{Thm25equ2}
	S(H^{\infty}_{v_{M,c}}(\CC))=\left\{(b_j)_{j\in\NN}\in\CC^{\NN}: \sup_{j\in\NN_{\ge 1}}\frac{M_{a_j+1}}{(\mu_{a_j+1})^{a_j+1}}\left(\sum_{a_j<l\le a_{j+1}}|b_l|^2\left(\frac{\mu_{a_j+1}}{c}\right)^{2l}\right)^{1/2}<+\infty\right\}.
\end{equation}
\end{theorem}

Concerning the solid core we get, by applying \cite[Theorem 2.4]{BonetLuskyTaskinen19}, the following characterization:

\begin{theorem}\label{Thm24}
Let $M\in\hyperlink{LCset}{\mathcal{LC}}$ be given. Assume that there exists a strictly increasing sequence (of integers) $(a_j)_{j\in\NN_{\ge 1}}$ and constants $b$ and $K$ with $K\ge b>2$ such that \eqref{Thm25equ} (the regularity condition $(b)$ holds true).

Then the solid core of $H^{\infty}_{v_{M,c}}(\CC)$ is given by
\begin{equation*}\label{Thm24equ1}
	s(H^{\infty}_{v_{M,c}}(\CC))=\left\{(b_j)_{j\in\NN}\in\CC^{\NN}: \sup_{j\in\NN_{\ge 1}}v_{M,c}(r_{a_j,c})\sum_{a_j<l\le a_{j+1}}|b_l|(r_{a_j,c})^{l}<+\infty\right\},
\end{equation*}
	or equivalently by
\begin{equation}\label{Thm24equ2}
	s(H^{\infty}_{v_{M,c}}(\CC))=\left\{(b_j)_{j\in\NN}\in\CC^{\NN}: \sup_{j\in\NN_{\ge 1}}\frac{M_{a_j+1}}{(\mu_{a_j+1})^{a_j+1}}\sum_{a_j<l\le a_{j+1}}|b_l|\left(\frac{\mu_{a_j+1}}{c}\right)^{l}<+\infty\right\}.
\end{equation}
\end{theorem}

Recall that, by the above computations, the (non)-existence of the Lusky numbers is not depending on the parameter $c>0$. Moreover, it is straightforward to see the following consequences:

\begin{corollary}\label{hullsequiv}
First, let $M\in\hyperlink{LCset}{\mathcal{LC}}$ be satisfying \eqref{Thm25equ} for a sequence $(a_j)_j$. Then
$$\forall\;c,d>0:\;\;\;S(H^{\infty}_{v_{M,c}}(\CC))\cong S(H^{\infty}_{v_{M,d}}(\CC)),\hspace{30pt}s(H^{\infty}_{v_{M,c}}(\CC))\cong s(H^{\infty}_{v_{M,d}}(\CC)),$$
i.e. all solid hulls and cores are isomorphic as sets.

Second, let $M,N\in\hyperlink{LCset}{\mathcal{LC}}$ be given with $M\hyperlink{approx}{\approx}N$ and such that $M$ or $N$ satisfy \eqref{Thm25equ} for some sequence $(a_j)_j$. Then we get the following isomorphisms as sets:
\begin{equation}\label{hullsequivequ}
\forall\;c,d>0:\;\;\;S(H^{\infty}_{v_{M,c}}(\CC))\cong S(H^{\infty}_{v_{N,d}}(\CC)),\hspace{30pt}s(H^{\infty}_{v_{M,c}}(\CC))\cong s(H^{\infty}_{v_{N,d}}(\CC)).
\end{equation}
\end{corollary}

\demo{Proof}
The isomorphism for the first part is realized via the mappings $(b_j)_j\mapsto((c/d)^jb_j)_j$ and $(b_j)_j\mapsto((d/c)^jb_j)_j$. For the second part, first the equivalence $M\hyperlink{approx}{\approx}N$ implies (see \eqref{Hinfinclusion}) that $\exists\;D\ge 1\;\forall\;c>0:$
$$S(H^{\infty}_{v_{N,c}}(\CC))\subseteq S(H^{\infty}_{v_{M,Dc}}(\CC))\subseteq S(H^{\infty}_{v_{N,D^2c}}(\CC)),\hspace{20pt}s(H^{\infty}_{v_{N,c}}(\CC))\subseteq s(H^{\infty}_{v_{M,Dc}}(\CC))\subseteq s(H^{\infty}_{v_{N,D^2c}}(\CC)).$$
Using this and the trivial inclusions $S(H^{\infty}_{v_{N,c}}(\CC))\subseteq S(H^{\infty}_{v_{N,c'}}(\CC))$, $s(H^{\infty}_{v_{N,c}}(\CC))\subseteq s(H^{\infty}_{v_{N,c'}}(\CC))$ for all $0<c\le c'$ (resp. for $M$) the conclusion follows by using the isomorphisms from the first part (for the sequence for which the numbers $(a_j)_j$ do exist).
\qed\enddemo

However, one shall note that for the different equivalence relation $\omega_M\hyperlink{sim}{\sim}\omega_N$ arising for weight functions in the ultradifferentiable setting, in general the analogous version of Corollary \ref{hullsequiv} will fail and also $H_{v_{M,c}}(\CC)\ncong H_{v_{N,d}}(\CC)$. But if $M$ does satisfy an additional technical assumption one is able to prove the following.

\begin{lemma}\label{hullsequiv1}
Let $M,N\in\hyperlink{LCset}{\mathcal{LC}}$ be given with $\omega_M\hyperlink{sim}{\sim}\omega_N$, such that $M$ or $N$ satisfies \eqref{Thm25equ} for some sequence $(a_j)_j$ and such that $M$ or $N$ has \hyperlink{mg}{$(\on{mg})$}. Then $M\hyperlink{approx}{\approx}N$ is valid and \eqref{hullsequivequ} holds true.
\end{lemma}

\demo{Proof}
First, by \cite[Proposition 3.6]{Komatsu73} property \hyperlink{mg}{$(\on{mg})$} for $M$ or $N$ is equivalent for $\omega_M$ or $\omega_N$ to have
\begin{equation}\label{omega6}
\exists\;H\ge 1\;\forall\;t\ge 0:\;\;\;2\omega(t)\le\omega(Ht)+H.
\end{equation}
For this condition, considered for arbitrary given weight functions, we refer to \cite[Corollary 16 $(3)$]{BonetMeiseMelikhov07}. Hence by equivalence \hyperlink{sim}{$\sim$} both $\omega_M$ and $\omega_N$ have \eqref{omega6} (e.g. see \cite[Lemma 3.2.2]{dissertation}) and the proof of \cite[Lemma 3.18 $(2)$]{testfunctioncharacterization} implies $M\hyperlink{approx}{\approx}N$. Thus the second part of Corollary \ref{hullsequiv} can be applied. More precisely, by iterating \eqref{omega6} one gets $\omega_M(t)\le C\omega_N(t)+C\le 2^n\omega_N(t)+C\le\omega_N(H_1t)+H_1$ for some $H_1\ge 1$ and all $t\ge 0$. This yields
$$\exists\;A,B\ge 1\;\forall\;c>0\;\forall\;t\ge 0:\;\;\;v_{M,Bc}(t)\le Av_{N,c}(t),\hspace{20pt}v_{N,Bc}(t)\le Av_{M,c}(t),$$
resp.
\begin{equation*}\label{hullsequiv1equ}
\exists\;B\ge 1\;\forall\;c>0:\;\;\;H^{\infty}_{v_{M,c}}(\CC)\subseteq H^{\infty}_{v_{N,Bc}}(\CC),\hspace{30pt}H^{\infty}_{v_{N,c}}(\CC)\subseteq H^{\infty}_{v_{M,Bc}}(\CC),
\end{equation*}
and this is precisely \eqref{Hinfinclusion}.
\qed\enddemo

To get more precise information on $S(H^{\infty}_{v_{M,c}}(\CC))$ and $s(H^{\infty}_{v_{M,c}}(\CC))$ via Theorems \ref{Thm25} and \ref{Thm24}, more precisely via \eqref{Thm25equ2} and \eqref{Thm24equ2}, we have to study the arising condition \eqref{Thm25equ}, i.e. the {\itshape regularity condition} $(b)$. Hence in our setting we are interested in asking:\vspace{6pt}

When $M\in\hyperlink{LCset}{\mathcal{LC}}$ is given, can we (always) find a strictly increasing sequence (of integers) $(a_j)_{j\in\NN_{\ge 1}}$, denoted by ''Lusky numbers'', such that there exist constants $K\ge b>2$ satisfying
 \begin{equation}\label{22}
\forall\;j\in\NN_{>0}:\;\;\;b\le\min\{A_M(a_j,a_{j+1}),B_M(a_j,a_{j+1})\}\le\max\{A_M(a_j,a_{j+1}),B_M(a_j,a_{j+1})\}\le K.
\end{equation}

Given such a sequence $(a_j)_j$, then it is clear that also each forward-shifted sequence $a^s:=(a_{j+s})_{j\ge 1}$, $s\in\NN_{>0}$ arbitrary, does satisfy \eqref{22} because $A_M(a^s_j,a^s_{j+1})=A_M(a_{j+s},a_{j+1+s})$, $B_M(a^s_j,a^s_{j+1})=B_M(a_{j+s},a_{j+1+s})$ for all $j\ge 1$. Related to this comment, on \cite[p. 598]{BonetTaskinen18} the following remark (translated into our weight sequence setting) has already been mentioned:

\begin{remark}\label{forwardshift}
Let $M\in\hyperlink{LCset}{\mathcal{LC}}$ be given such that there exists a strictly increasing sequence (of integers) $(a_j)_{j\in\NN_{\ge 1}}$ such that \eqref{Thm25equ} holds true. Then we can replace in \eqref{Thm25equ1} $v_{M,c}(r_{a_j})$ by $v_{M,c}(r_{a_{j+1}})$ and $r_{a_j}$ by $r_{a_{j+1}}$ resp. in \eqref{Thm25equ2} $\frac{M_{a_j+1}}{(\mu_{a_j+1})^{a_j+1}}$ by $\frac{M_{a_{j+1}+1}}{(\mu_{a_{j+1}+1})^{a_{j+1}+1}}$ and $\mu_{a_j+1}$ by $\mu_{a_{j+1}+1}$ (and similarly in Theorem \ref{Thm24}).

Note that in the representations of $S(H^{\infty}_{v_{M,c}}(\CC))$ and $s(H^{\infty}_{v_{M,c}}(\CC))$ equivalently we can start the summation at any $a_{j_0}$, $j_0>1$.
\end{remark}

Moreover we see the following observations:

\begin{remark}\label{firstremark}
\begin{itemize}
\item[$(i)$] It is immediate that we never can choose $(a_j)_j$ to be constant and it is also impossible to have $a_{j+1}=a_j+1$, since in this case $B_M(a_j,a_{j+1})=1$ follows. Thus a necessary growth assumption is
    $$a_{j+1}\ge a_j+2,\;\;\;\forall\;j\ge 1.$$

\item[$(ii)$] If $M\in\hyperlink{LCset}{\mathcal{LC}}$ satisfies \eqref{22} with $(a_j)_j$, then $M^s$ as well for any $s>1$ because clearly $A_{M^s}(a_j,a_{j+1})=(A_M(a_j,a_{j+1}))^s$, $B_{M^s}(a_j,a_{j+1})=(B_M(a_j,a_{j+1}))^s$. (For $0<s<1$, in general this will be not valid anymore since for small $s>0$ the first estimate can fail).

    More generally, if $M,N\in\hyperlink{LCset}{\mathcal{LC}}$ both are satisfying \eqref{22} with the same sequence $(a_j)_j$, then $(a_j)_j$ can be used for $M\cdot N:=(M_pN_p)_p$ as well: The quotients of this sequence are given by $\mu_j\nu_j$ and so $A_{M\cdot N}(a_j,a_{j+1})=A_M(a_j,a_{j+1})A_N(a_j,a_{j+1})$, $B_{M\cdot N}(a_j,a_{j+1})=B_M(a_j,a_{j+1})B_N(a_j,a_{j+1})$.

\item[$(iii)$] For the sequence $Q:=M/N:=(M_p/N_p)_p$ we can only show the upper bound: One has $\rho_j=\mu_j/\nu_j$ and so $A_Q(a_j,a_{j+1})=A_M(a_j,a_{j+1})\frac{\nu_{a_j+1}\nu_{a_j+2}\cdots\nu_{a_{j+1}}}{(\nu_{a_j+1}+1)^{a_{j+1}-a_j}}\le A_M(a_j,a_{j+1})\frac{\nu_{a_j+1}\nu_{a_j+2}\cdots\nu_{a_{j+1}}}{(\nu_{a_j+1})^{a_{j+1}-a_j}}=A_M(a_j,a_{j+1})B_N(a_j,a_{j+1})$, similarly $B_Q(a_j,a_{j+1})\le B_M(a_j,a_{j+1})A_N(a_j,a_{j+1})$.
\end{itemize}
\end{remark}

\subsection{On the $r$-interpolating sequence and ramified weights}\label{rinterpolatingsect}
Given $M\in\hyperlink{LCset}{\mathcal{LC}}$ satisfying \eqref{22} with $(a_j)_j$, then $(ii)$ in Remark \ref{firstremark} provides a first method to construct new sequences still satisfying \eqref{22} (with the same Lusky numbers $(a_j)_j$). Another possibility is the following approach:

In \cite[Lemma 2.3]{Schmetsvaldivia00} for given $M\in\hyperlink{LCset}{\mathcal{LC}}$ and $r\in\NN_{\ge 1}$ the so-called $r$-interpolating sequence $P^{M,r}$ has been introduced as follows, see also \cite[Sect. 2.5]{mixedramisurj}:
\begin{equation}\label{interpolatingsequ}
P^{M,r}_{rk+j}:=((M_k)^{r-j}(M_{k+1})^j)^{1/r},\hspace{20pt}\forall\;k\in\NN,\;\forall\;j\in\{0,\dots,r\}.
\end{equation}
We have $P^{M,r}_{rj}=M_j$ for all $j\in\NN$ (i.e. we get $P^{M,1}\equiv M$) and by denoting $\pi^{M,r}_k:=\frac{P^{M,r}_k}{P^{M,r}_{k-1}}$ (with $\pi^{M,r}_0:=1$) we see
\begin{equation}\label{interpolatingsequ0}
\forall\;k\in\NN\;\forall\;j\in\{1,\dots, r\}:\;\;\;\pi^{M,r}_{rk+j}=\frac{P^{M,r}_{rk+j}}{P^{M,r}_{rk+j-1}}=\left(\frac{M_k^{r-j}M_{k+1}^j}{M_k^{r-j+1}M_{k+1}^{j-1}}\right)^{1/r}=\left(\frac{M_{k+1}}{M_k}\right)^{1/r}=(\mu_{k+1})^{1/r}.
\end{equation}
Hence $M\in\hyperlink{LCset}{\mathcal{LC}}$ if and only if $P^{M,r}\in\hyperlink{LCset}{\mathcal{LC}}$ and by using \eqref{interpolatingsequ0} we prove the following observation.

\begin{lemma}\label{rinterpolatinglemma}
Let $M\in\hyperlink{LCset}{\mathcal{LC}}$ be given and $r\in\NN_{\ge 1}$. Assume that there exists a sequence of integers $(a_j)_j$ satisfying \eqref{22} (for $M$). Then the $r$-interpolating sequence $P^{M,r}$ does satisfy \eqref{22} for the choice $(ra_j)_j$, i.e. stretching the Lusky numbers by the factor $r$.
\end{lemma}

\demo{Proof}
For simplicity we write in the proof $\pi_k$ instead of $\pi^{M,r}_k$.

First, we have $(\pi_{ra_{j+1}+1})^{r(a_{j+1}-a_j)}=(\mu_{a_{j+1}+1})^{a_{j+1}-a_j}$.

Second, we decompose the product $\pi_{ra_j+1}\pi_{ra_j+2}\cdots\pi_{ra_{j+1}}$ into $a_{j+1}-a_j$ many factors of length $r$ as follows: Since $\pi_{ra_j+i}=(\mu_{a_j+1})^{1/r}$ for all $1\le i\le r$, we get $\pi_{ra_j+1}\cdots\pi_{ra_j+r}=\mu_{a_j+1}$. Similarly $\pi_{r(a_j+1)+i}=(\mu_{a_j+2})^{1/r}$ for all $1\le i\le r$, which yields $\pi_{ra_j+r+1}\cdots\pi_{ra_j+2r}=\mu_{a_j+2}$. Finally, in the last block we have $\pi_{r(a_{j+1}-1)+i}=(\mu_{a_{j+1}})^{1/r}$ for all $1\le i\le r$, which yields $\pi_{ra_{j+1}-r+1}\cdot\pi_{ra_{j+1}-r+2}\cdots\pi_{ra_{j+1}}=\mu_{a_{j+1}}$.

Altogether, by \eqref{Ainteresting} we have shown $$A_M(a_j,a_{j+1})=\frac{(\pi_{ra_{j+1}+1})^{r(a_{j+1}-a_j)}}{\pi_{ra_j+1}\pi_{ra_j+2}\cdots\pi_{ra_{j+1}}}=A_{P^{M,r}}(ra_j,ra_{j+1}).$$

Similarly we have $(\mu_{a_j+1})^{a_{j+1}-a_j}=(\pi_{ra_j+1})^{r(a_{j+1}-a_j)}$ and since the product arising in the numerator of $B_M(a_j,a_{j+1})$ is precisely the same as in the denominator of $A_M(a_j,a_{j+1})$ we have $B_M(a_j,a_{j+1})=B_{P^{M,r}}(ra_j,ra_{j+1})$ as well.
\qed\enddemo

In addition, the $r$-interpolating sequence can be used to get some information on the inclusion of spaces of weighted entire functions w.r.t. relation $\omega_N\hyperlink{ompreceq}{\preceq}\omega_M$:

Let $M,N\in\hyperlink{LCset}{\mathcal{LC}}$ such that $\omega_N\hyperlink{ompreceq}{\preceq}\omega_M$ holds true, then $\omega_M(t)\le r\omega_N(t)+r$ for all $t\ge 0$ and some $r\ge 1$ and w.l.o.g. we can take $r\in\NN_{\ge 1}$. One has for all $t\ge 0$ that
\begin{align*}
r\omega_N(t)&=r\sup_{j\in\NN}\log\left(\frac{t^j}{N_j}\right)=\sup_{j\in\NN}\log\left(\frac{t^{rj}}{(N_j)^r}\right)=\omega_{(N)^r}(t^r)\le\sup_{j\in\NN}\log\left(\frac{t^{rj}}{N_j}\right)=\sup_{j\in\NN}\log\left(\frac{t^{rj}}{P^{N,r}_{rj}}\right)
\\&
\le\omega_{P^{N,r}}(t),
\end{align*}
which implies
\begin{equation}\label{Pinterpolatingequ1}
\forall\;c>0:\;\;\;H^{\infty}_{v_{M,c}}(\CC)\subseteq H^{\infty}_{v_{P^{N,r},c}}(\CC).
\end{equation}
Consequently, $\omega_M\hyperlink{sim}{\sim}\omega_N$ yields
\begin{equation}\label{Pinterpolatingequ2}
\exists\;r\in\NN_{\ge 1}\;\forall\;c>0:\;\;\;H^{\infty}_{v_{M,c}}(\CC)\subseteq H^{\infty}_{v_{P^{N,r},c}}(\CC),\;\;\;H^{\infty}_{v_{N,c}}(\CC)\subseteq H^{\infty}_{v_{P^{M,r},c}}(\CC).
\end{equation}

Moreover, the sequence $P^{M,r}$ can be used to see how ramification of the complex variable is translated into the weight sequence setting.

\begin{lemma}\label{ramification}
Let $M\in\hyperlink{LCset}{\mathcal{LC}}$ be given, $r\in\NN_{\ge 1}$ and $P^{M,r}$ the according $r$-interpolating sequence. Then for all $t\ge 0$ we get
$$\omega_M(t^r)=r^2\omega_{P^{M,r}}(t).$$
\end{lemma}

\demo{Proof}
We use the following integral representation formula for $\omega_M$, see \cite[1.8.III]{mandelbrojtbook}:
$$\omega_M(t)=\int_0^{t}\frac{\Sigma_M(u)}{u}du=\int_{\mu_1}^{t}\frac{\Sigma_M(u)}{u}du,$$
with $\Sigma_M(t):=|\{p\in\NN_{>0}:\;\mu_p\le t\}|=\max\{p\in\NN_{>0}:\;\mu_p\le t\}$. Then by \eqref{interpolatingsequ0} we get $$\Sigma_M(t^r)=r\Sigma_{P^{M,r}}(t)$$
for all $t\ge 0$: If $\mu_{k+1}\le t<\mu_{k+2}$ for some $k\in\NN$, then $\Sigma_M(t)=k+1$ and $\pi^{M,r}_{rk+j}=(\mu_{k+1})^{1/r}\le t^{1/r}<(\mu_{k+2})^{1/r}=\pi^{M,r}_{r(k+1)+j}$ for all $1\le j\le r$ does precisely give $\Sigma_{P^{M,r}}(t^{1/r})=r(k+1)$.

Note that $(\mu_1)^{1/r}=\pi^{M,r}_j$ for all $1\le j\le r$, thus $\Sigma_{P^{M,r}}(t)=0$ for $0\le t<(\mu_1)^{1/r}$ and so precisely for all $t$ satisfying $0\le t^r<\mu_1$, i.e. for all $t$ satisfying $\Sigma_M(t^r)=0$.

Using this we can calculate as follows:
\begin{align*}
\omega_M(t^r)&=\int_{\mu_1}^{t^r}\frac{\Sigma_M(x)}{x}dx=\int_{(\mu_1)^{1/r}}^{t}\frac{\Sigma_M(y^r)}{y^r}ry^{r-1}dy=r\int_{(\mu_1)^{1/r}}^{t}\frac{\Sigma_M(y^r)}{y}dy
\\&
=r^2\int_{(\mu_1)^{1/r}}^t\frac{\Sigma_{P^{M,r}}(y)}{y}dy=r^2\int_{\pi^{M,r}_1}^t\frac{\Sigma_{P^{M,r}}(y)}{y}dy=r^2\omega_{P^{M,r}}(t).
\end{align*}
\qed\enddemo

For any given weight function $v$ and $r>0$ we set
\begin{equation}\label{ramifiedweight}
v^r(t):=v(t^r),\hspace{30pt}{}^rv(t):=(v(t))^r,
\end{equation}
so $v^1\equiv v\equiv{}^1v$. It is immediate that each $v^r$, ${}^rv$ is again a weight function because $\lim_{t\rightarrow+\infty}t^kv^r(t)=\lim_{t\rightarrow+\infty}t^kv(t^r)=\lim_{s\rightarrow+\infty}s^{k/r}v(s)=0$ for all $k\ge 0$ and similarly for ${}^rv$. We set
$$H^{\infty}_{v^r}(\CC):=\{f\in H(\CC): \|f\|_{v^r}:=\sup_{z\in\CC}|f(z)|v^r(|z|)<+\infty\},$$
and similarly for ${}^rv$. When $r\in\NN_{>0}$ and $v\equiv v_{M,c}$, then the $r$-interpolating sequence $P^{M,r}$ can be used to determine the solid hull and solid core of $H^{\infty}_{w_{M,c}^r}(\CC)$, where we set $w^r_{M,c}:t\mapsto\exp(-r^{-2}\omega_{M}((ct)^r))$.

\begin{theorem}\label{Pramification}
Let $M\in\hyperlink{LCset}{\mathcal{LC}}$ be given and $r\in\NN_{\ge 1}$. Assume that there exists a sequence of integers $(a_j)_j$ satisfying \eqref{22} (the regularity condition $(b)$), then
$$S(H^{\infty}_{w^r_{M,c}}(\CC))=\left\{(b_j)_{j\in\NN}\in\CC^{\NN}: \sup_{j\in\NN_{\ge 1}}\frac{M_{a_j+1}}{(\mu_{a_j+1})^{a_j+1}}\left(\sum_{ra_j<l\le ra_{j+1}}|b_l|^2\left(\frac{(\mu_{a_j+1})^{1/r}}{c}\right)^{2l}\right)^{1/2}<+\infty\right\},$$
and
$$s(H^{\infty}_{w^r_{M,c}}(\CC))=\left\{(b_j)_{j\in\NN}\in\CC^{\NN}: \sup_{j\in\NN_{\ge 1}}\frac{M_{a_j+1}}{(\mu_{a_j+1})^{a_j+1}}\sum_{ra_j<l\le ra_{j+1}}|b_l|\left(\frac{(\mu_{a_j+1})^{1/r}}{c}\right)^{l}<+\infty\right\}.$$
\end{theorem}

\demo{Proof}
First, Lemma \ref{ramification} does imply $H^{\infty}_{w^r_{M,c}}(\CC)=H^{\infty}_{v_{P^{M,r}},c}(\CC)$. Lemma \ref{rinterpolatinglemma}
yields that $(ra_j)_j$ satisfies \eqref{22} for $P^{M,r}$.

Now write $P$ instead of $P^{M,r}$, then \eqref{Thm25equ2} yields
 $$S(H^{\infty}_{w^r_{M,c}}(\CC))=\left\{(b_j)_{j\in\NN}\in\CC^{\NN}: \sup_{j\in\NN_{\ge 1}}\frac{P_{ra_j+1}}{(\pi_{ra_j+1})^{ra_j+1}}\left(\sum_{ra_j<l\le ra_{j+1}}|b_l|^2\left(\frac{\pi_{ra_j+1}}{c}\right)^{2l}\right)^{1/2}<+\infty\right\},$$
and \eqref{Thm24equ2} yields
$$s(H^{\infty}_{w^r_{M,c}}(\CC))=\left\{(b_j)_{j\in\NN}\in\CC^{\NN}: \sup_{j\in\NN_{\ge 1}}\frac{P_{ra_j+1}}{(\pi_{ra_j+1})^{ra_j+1}}\sum_{ra_j<l\le ra_{j+1}}|b_l|\left(\frac{\pi_{ra_j+1}}{c}\right)^{l}<+\infty\right\}.$$
Finally, \eqref{interpolatingsequ} and \eqref{interpolatingsequ0} imply that
$$\frac{P_{ra_j+1}}{(\pi_{ra_j+1})^{ra_j+1}}=\frac{((M_{a_j})^{r-1}M_{a_j+1})^{1/r}}{(\mu_{a_j+1})^{a_j+1/r}}=\frac{M_{a_j}(\mu_{a_j+1})^{1/r}}{(\mu_{a_j+1})^{a_j+1/r}}=\frac{M_{a_j}}{(\mu_{a_j+1})^{a_j}}=\frac{M_{a_j+1}}{(\mu_{a_j+1})^{a_j+1}},$$
and $\pi_{ra_j+1}=(\mu_{a_j+1})^{1/r}$.
\qed\enddemo

For the sake of completeness, we finish this section with the following comments on arbitrary weight functions $v$ and ramification parameters $r>1$:

\begin{itemize}
\item[$(i)$] For any $p>0$, when $t_p$ is denoting the global maximum point of $t\mapsto t^pv(t)$, i.e. $t_p^pv(t_p)\ge t^pv(t)$ for all $t\ge 0$, we get $t_p^{rp}(v(t_p))^r\ge t^{rp}(v(t))^r$ for all $t\ge 0$.

    So with $s:=t^r$ and $s_p:=t_p^r$ one has $s_p^p\cdot{}^rv((s_p)^{1/r})\ge s^p\cdot{}^rv(s^{1/r})\Leftrightarrow s_p^p\cdot{}^rv^{1/r}(s_p)\ge s^p\cdot{}^rv^{1/r}(s)$ for all $s\ge 0$, i.e. $s_p$ is the maximum point of $s\mapsto s^p\cdot{}^rv^{1/r}(s)$.

\item[$(ii)$] Thus, for the expressions under consideration in the regularity condition $(b)$ for $v$, see \cite[$(2.1)$]{BonetTaskinen18}, we get for any $0<m<n$:
$$A(m,n):=\left(\frac{t_m}{t_n}\right)^m\frac{v(t_m)}{v(t_n)}=\left(\frac{s_m}{s_n}\right)^{m/r}\left(\frac{{}^rv^{1/r}(s_m)}{{}^rv^{1/r}(s_n)}\right)^{1/r},$$
$$B(m,n):=\left(\frac{t_n}{t_m}\right)^n\frac{v(t_n)}{v(t_m)}=\left(\frac{s_n}{s_m}\right)^{n/r}\left(\frac{{}^rv^{1/r}(s_n)}{{}^rv^{1/r}(s_m)}\right)^{1/r}.$$
Thus, if $r>1$ and $v$ does satisfy the regularity condition $(b)$, then each weight ${}^rv^{1/r}$ as well with the same sequence of Lusky numbers as for $v$. (If $0<r<1$, then in general the estimate from above will fail.)
\end{itemize}

\subsection{First (counter)-example}\label{firstcounter}
Next let us see that not each sequence $M\in\hyperlink{LCset}{\mathcal{LC}}$ does automatically have \eqref{22}, i.e. not for each weight sequence there do exist the ''Lusky numbers''.

\begin{lemma}\label{whatcanhappen}
There does exist $M\in\hyperlink{LCset}{\mathcal{LC}}$ such that there does not exist a sequence of integers $(a_j)_j$ satisfying \eqref{22}.
\end{lemma}

\demo{Proof}
We define $M$ in terms of the sequence of quotients $(\mu_p)_{p\in\NN}$, i.e. $M_p=\prod_{i=1}^p\mu_i$, $p\in\NN$, with $\mu_0:=1$ and $\mu_p\rightarrow+\infty$ as $p\rightarrow+\infty$.

We set $\mu_0:=1$ and let $(b_j)_{j\in\NN_{>0}}$ be an arbitrary strictly increasing sequence in $\NN$ such that $b_1=1$. Then put
$$\mu_p:=c_j,\hspace{20pt}\;b_j\le p<b_{j+1},$$
with $(c_j)_{j\in\NN_{>0}}$ an arbitrary sequence of positive real numbers satisfying $1\le c_1$ (normalization), $c_j<c_{j+1}$ and $\lim_{j\rightarrow+\infty}\frac{c_{j+1}}{c_j}=+\infty$.\vspace{6pt}

Now take a given strictly increasing sequence $(a_j)_j$ (of integers) with $a_1\ge 1$, and we analyze $A_M(a_j,a_{j+1})$. Given $a_j$, $j\ge 1$, we have $b_{l_j}\le a_j+1<b_{l_{j+1}}$ for some $l_j\in\NN$.

We distinguish now between two cases: If also $b_{l_j}\le a_{j+1}+1<b_{l_j+1}$, then $A_M(a_j,a_{j+1})=1$ since in the numerator and in the denominator we have the same product $(c_{l_j})^{a_{j+1}-a_j}$.

If $a_{j+1}+1\ge b_{l_j+1}$, then $b_{l_j+k_j}\le a_{j+1}+1<b_{l_j+k_j+1}$ for some $k_j\in\NN_{\ge 1}$ and we set $d_j:=\frac{c_{l_j+1}}{c_{l_j}}$. In this situation either $b_{l_j+k_j}\le a_{j+1}$ or $b_{l_j+k_j}=a_{j+1}+1$ and in both cases we have that
\begin{align*}
A_M(a_j,a_{j+1})\ge\left(\frac{c_{l_j+k_j}}{c_{l_j}}\right)^{b_{l_j+1}-a_j-1}\ge\frac{c_{l_j+1}}{c_{l_j}}=d_j.
\end{align*}


Since $d_j\rightarrow+\infty$ as $j\rightarrow+\infty$ we see that for any choice $(a_j)_j$ at least one estimate in \eqref{22} has to fail.
\qed\enddemo

Choosing the sequences $(b_j)_j$ and $(c_j)_j$ in a more precise (convenient) way we see that such constructed sequences can satisfy several known growth and regularity properties used in the theory of ultradifferentiable (and ultraholomorphic) functions.

\begin{corollary}\label{whatcanhappencor}
There does exist $M\in\hyperlink{LCset}{\mathcal{LC}}$ satisfying \hyperlink{beta1}{$(\beta_1)$} and \hyperlink{dc}{$(\on{dc})$}, i.e. $M$ is strongly non-quasianalytic and has derivation closedness, but such that there does not exist a sequence of integers $(a_j)_j$ satisfying \eqref{22}.
\end{corollary}

\demo{Proof}
Let $Q\in\NN_{\ge 2}$ and $D>1$ and then take
$$b_1:=1,\;\;\;b_{j+1}:=Qb_j,\;\hspace{30pt}c_j:=D^{b_j},\;\;\;j\in\NN_{\ge 1}.$$
Hence $b_{j+1}=Q^j$ for $j\in\NN_{\ge 1}$ and $c_1=D$, $c_j=D^{Q^{j-1}}$, $j\in\NN_{\ge 2}$, and we see that $\frac{c_{j+1}}{c_j}=D^{b_{j+1}-b_j}=D^{Q^{j-1}(Q-1)}\rightarrow+\infty$ as $j\rightarrow+\infty$.

Moreover $M$ does have \hyperlink{beta1}{$(\beta_1)$} (hence \hyperlink{beta3}{$(\beta_3)$} too) because $\frac{\mu_{Qp}}{\mu_p}=\frac{c_{j+1}}{c_j}$ for all $p\in\NN$ with $b_j\le p<b_{j+1}$ and \hyperlink{dc}{$(\text{dc})$} follows because $\mu_p\le D^p$ for all $p\in\NN$ by definition.

However, $M$ does not have \hyperlink{mg}{$(\text{mg})$} since for this property it is required to have $\sup_{p\in\NN}\frac{\mu_{2p}}{\mu_p}<\infty$ (e.g. see \cite[Lemma 2.2]{whitneyextensionweightmatrix}) and this property is obviously violated. But in Lemma \ref{whatcanhappen1} below, by using the techniques developed in Section \ref{Alternativesequencedelta}, we provide an analogous example even satisfying \hyperlink{mg}{$(\text{mg})$}.
\qed\enddemo

\subsection{First example - $q$-Gevrey sequences}\label{qGevrey}
The aim is now to give a more concrete example for the representations obtained in Theorems \ref{Thm25} and \ref{Thm24}. We study the (family of) sequences
$$M^q=(M^q_p)_{p\in\NN},\hspace{30pt}M^q_p:=q^{p^2},\;\;q>1,$$
hence $\mu^q_p:=\frac{M^q_p}{M^q_{p-1}}=q^{2p-1}$ for all $p\in\NN_{\ge 1}$ (and set $\mu_0:=1$). Each $M^q$ does have \hyperlink{beta1}{$(\beta_1)$} and \hyperlink{dc}{$(\on{dc})$}, but none of them has \hyperlink{mg}{$(\on{mg})$}.\vspace{6pt}

The goal is to prove the following result:

\begin{proposition}\label{qGevreyprop}
	For $M^q$ property \eqref{22} is satisfied for any sequence (of integers) $(a_j)_{j\in\NN_{\ge 1}}$ satisfying
	$$\sup_{j\in\NN_{\ge 1}}a_{j+1}-a_j<+\infty,\hspace{30pt}(a_{j+1}-a_j)(a_{j+1}-a_j-1)>\frac{\log(2)}{\log(q)},\;j\in\NN_{>0}.$$
\end{proposition}
In particular, when given $q>1$, then with $c$ chosen large enough to guarantee $q>2^{1/(c(c-1))}$ (and $c\in\NN_{\ge 2}$) we can take $a_{j+1}=a_j+c$, $a_1:=1$, i.e.
\begin{equation}\label{qGevreychoice}
a_j=c(j-1)+1.
\end{equation}

\demo{Proof}
From now on we fix the parameter $q$ and for simplicity we write $\mu_p$ instead of $\mu^q_p$. First we study the expression $A_M(a_j,a_{j+1})$ in \eqref{Ainteresting}.

One has
$$(\mu_{a_{j+1}+1})^{a_{j+1}-a_j}=q^{(2(a_{j+1}+1)-1)(a_{j+1}-a_j)}$$
and
$$\mu_{a_j+1}\cdots\mu_{a_{j+1}}=q^{2(a_j+1)-1+2(a_j+2)-1+\dots+2a_{j+1}-1}=q^{2(a_j+1)+\dots+2a_{j+1}-(a_{j+1}-a_j)}.$$
More precisely, the first part of the argument in this exponent yields $2(a_j+1)+2(a_j+2)+\dots+2a_{j+1}=2a_j(a_{j+1}-a_j)+2\sum_{l=1}^{a_{j+1}-a_j}l=2a_j(a_{j+1}-a_j)+2\frac{(a_{j+1}-a_j)(a_{j+1}-a_j+1)}{2}$.

Altogether we have for the argument arising in the exponent of $A_M(a_j,a_{j+1})$:
\begin{align*}
&(2(a_{j+1}+1)-1)(a_{j+1}-a_j)-2a_j(a_{j+1}-a_j)-(a_{j+1}-a_j)(a_{j+1}-a_j+1)+(a_{j+1}-a_j)
\\&
=2(a_{j+1}+1)(a_{j+1}-a_j)-2a_j(a_{j+1}-a_j)-(a_{j+1}-a_j)(a_{j+1}-a_j+1)
\\&
=2(a_{j+1}-a_j)(a_{j+1}-a_j+1)-(a_{j+1}-a_j)(a_{j+1}-a_j+1)=(a_{j+1}-a_j)(a_{j+1}-a_j+1).
\end{align*}
Thus $A_M(a_j,a_{j+1})=q^{(a_{j+1}-a_j)(a_{j+1}-a_j+1)}$ and so the upper estimate of \eqref{22} holds if and only if $\sup_{j\in\NN_{\ge 1}}a_{j+1}-a_j<+\infty$. The lower estimate of \eqref{22} holds if and only if $(a_{j+1}-a_j)(a_{j+1}-a_j+1)>\log(2)/\log(q)$ for all $j\in\NN_{\ge 1}$.\vspace{6pt}

Next we study $B_M(a_j,a_{j+1})$ in \eqref{Binteresting}: Similarly as above we get
$$\mu_{a_j+1}\mu_{a_j+2}\cdots\mu_{a_{j+1}}=q^{2(a_j+1)-1+2(a_j+2)-1+\dots+2a_{j+1}-1}=q^{2(a_j+1)+\dots+2a_{j+1}-(a_{j+1}-a_j)}$$
and $(\mu_{a_j+1})^{a_{j+1}-a_j}=q^{(2(a_j+1)-1)(a_{j+1}-a_j)}$. Hence the argument arising in the exponent of $B_M(a_j,a_{j+1})$ is given by
\begin{align*}
&2(a_j+1)+\dots+2a_{j+1}-(a_{j+1}-a_j)-(2(a_j+1)-1)(a_{j+1}-a_j)
\\&
=2a_j(a_{j+1}-a_j)+(a_{j+1}-a_j)(a_{j+1}-a_j+1)-(a_{j+1}-a_j)-(2(a_j+1)-1)(a_{j+1}-a_j)
\\&
=2a_j(a_{j+1}-a_j)+(a_{j+1}-a_j)(a_{j+1}-a_j+1)-2(a_j+1)(a_{j+1}-a_j)
\\&
=(a_{j+1}-a_j)(2a_j+a_{j+1}-a_j+1-2a_j-2)=(a_{j+1}-a_j)(a_{j+1}-a_j-1).
\end{align*}
Thus $B_M(a_j,a_{j+1})=q^{(a_{j+1}-a_j)(a_{j+1}-a_j-1)}$ and so the upper estimate of \eqref{22} holds again if and only if $\sup_{j\in\NN_{\ge 1}}a_{j+1}-a_j<+\infty$. The lower estimate of \eqref{22} holds if and only if $(a_{j+1}-a_j)(a_{j+1}-a_j-1)>\log(2)/\log(q)$ for all $j\in\NN_{\ge 1}$.
\qed\enddemo

For any $q\ge 2$ we can choose $c=2$ in Proposition \ref{qGevreyprop}, so we obtain via \eqref{qGevreychoice}
$$a_j=2j-1,\hspace{30pt}\frac{M_{a_j+1}}{(\mu_{a_j+1})^{a_j+1}}=q^{-a_j^2-a_j}=q^{-2j(2j-1)},\hspace{30pt}\mu_{a_j+1}=q^{4j-1}.$$
Thus Theorems \ref{Thm25} and \ref{Thm24} give the following characterizations:

\begin{corollary}\label{Thm25qgevrey}
	Let $M^q$ be given with $q\ge 2$. Then the solid hull of $H^{\infty}_{v_{M^q,c}}$ is given by the set
	$$S(H^{\infty}_{v_{M^q,c}}(\CC))=\left\{(b_j)_{j\in\NN}\in\CC^{\NN}: \sup_{j\in\NN_{\ge 1}}q^{-2j(2j-1)}\left(\sum_{2j\le l\le 2j+1}|b_l|^2\left(\frac{q^{4j-1}}{c}\right)^{2l}\right)^{1/2}<+\infty\right\},$$
	or equivalently in a more compact form:
	$$S(H^{\infty}_{v_{M^q,c}}(\CC))=\left\{(b_j)_{j\in\NN}\in\CC^{\NN}: \exists\;D\ge 1\;\forall\;j\in\NN_{\ge 1}:\;\;\;|b_{2j}|^2+|b_{2j+1}|^2\frac{q^{2(4j-1)}}{c^2}\le D^2c^{4j}q^{-8j^2}\right\}.$$
Moreover we get
$$s(H^{\infty}_{v_{M^q,c}}(\CC))=\left\{(b_j)_{j\in\NN}\in\CC^{\NN}: \exists\;D\ge 1\;\forall\;j\in\NN_{\ge 1}:\;\;\;|b_{2j}|+|b_{2j+1}|\frac{q^{4j-1}}{c}\le Dc^{2j}q^{-4j^2}\right\}.$$
\end{corollary}

\demo{Proof}
Concerning the solid hull, the first identity is immediate by the possible choice of $(a_j)_j$. For each sequence $(b_j)_j$ contained in the first set it is equivalent that there exists some $D\ge 1$ such that for all $j\ge 1$ we have
\begin{align*}
&|b_{2j}|^2\left(\frac{q^{4j-1}}{c}\right)^{4j}+|b_{2j+1}|^2\left(\frac{q^{4j-1}}{c}\right)^{4j+2}\le D^2q^{4j(2j-1)}\Longleftrightarrow|b_{2j}|^2+|b_{2j+1}|^2\frac{q^{8j-2}}{c^2}\le D^2c^{4j}q^{-8j^2}.
\end{align*}
The solid core follows analogously.
\qed\enddemo

By using the explicit representations obtained for the solid hull and solid core we can prove:

\begin{corollary}\label{Thm25qgevrey1}
	Let $M^q$ be given with $q\ge 2$. Then
\begin{equation}\label{Gevreyhullcore}
\forall\;c>0:\;\;\;S(H^{\infty}_{v_{M^q,c}}(\CC))=H^{\infty}_{v_{M^q,c}}(\CC)=s(H^{\infty}_{v_{M^q,c}}(\CC)),
\end{equation}
so the space $H^{\infty}_{v_{M^q,c}}(\CC)$ is solid.
\end{corollary}

\demo{Proof}
\eqref{Gevreyhullcore} follows by using the representations from Corollary \ref{Thm25qgevrey}. For this we point out that for all $a,b,c\ge 0$ the inequality $a+b\le c$ implies $a^2+b^2\le c^2$. Conversely, $a^2+b^2\ge\frac{1}{2}(a+b)^2\Leftrightarrow(a-b)^2\ge 0$ is valid and so $a^2+b^2\le c^2$ implies $a+b\le\sqrt{2}c$.
\qed\enddemo

We close this section by summarizing some facts for the weight(s) $v_{M^q,c}$.

\begin{remark}\label{logsqureremark}
\begin{itemize}

\item[$(i)$] By \cite[Lemma 5.7]{compositionpaper} and \cite[Sect. 5.5]{whitneyextensionweightmatrix} we get that each $\omega_{M^q}$ is equivalent (w.r.t. \hyperlink{sim}{$\sim$}) to the weight $t\mapsto\max\{0,(\log(t))^2\}$, alternatively also to $t\mapsto(\log(1+t))^2$.

\item[$(ii)$] For this abstractly given weight (which is violating \eqref{omega6}) the solid hull and solid core has not been computed before in terms of the Lusky numbers. However, the weight $v(t):=\exp(-\max\{0,(\log(t))^2\})$ resp. $v(t):=\exp(-(\log(1+t))^2)$ has been studied in the literature, see \cite[Thm. 2.4, Cor.]{Lusky00}, \cite[Sect. 5]{BonetLuskyTaskinen19Bergman} and \cite[Sect. 5]{BonetLuskyTaskinen19}.

\item[$(iii)$] More precisely, in \cite[Sect. 5, Example]{BonetLuskyTaskinen19} by using \cite[Thm. 2.4, Cor.]{Lusky00} it has been shown that $H^{\infty}_{v}(\CC)$ is solid which should be compared with Corollary \ref{Thm25qgevrey1}.

\item[$(iv)$] More generally, in \cite[Theorem 5.2]{BonetLuskyTaskinen19Bergman} it has been shown that for any weight $v$ the space $H^{\infty}_{v}(\CC)$ is solid if and only if $\sup_{j\ge 1}a_{j+1}-a_j<+\infty$ (see also $(ii)$ in Lemma \ref{alternativedeltaconsequ} below).

\item[$(v)$] This last comment should be compared with the admissible choice for the Lusky numbers in \eqref{qGevreychoice} obtained before which shows that our approach is consistent with the known results.
\end{itemize}
\end{remark}

\section{Alternative representations for the regularity condition $(b)$}\label{Alternativesequencedelta}

In this section we derive alternative useful representations for the expressions $A_M(a_j,a_{j+1})$ and $B_M(a_j,a_{j+1})$. As we see this method is convenient to get more information on the existence and (possible) growth behavior of the Lusky numbers $a_j$ in the weight sequence setting resp. how the growth of the Lusky numbers is related to the growth of the sequence $M$.

This approach has been inspired by the construction of (counter)-examples (see \cite[Sect. 2.2.5]{dissertationjimenez}), for similar computations see also e.g. \cite[Prop. 3.3]{logconvexnonproximate}.\vspace{6pt}

For any given $M\in\hyperlink{LCset}{\mathcal{LC}}$ (recall $\mu_0:=1$) we put
\begin{equation*}\label{deltachoice}
\delta_p:=\log(\mu_p)-\log(\mu_{p-1})=\log(\mu_p/\mu_{p-1}),\;\;\;p\ge 1,
\end{equation*}
hence
\begin{equation}\label{muchoice}
\mu_p=\exp\left(\sum_{j=1}^p\delta_j\right),\;\;\;p\ge 1.
\end{equation}
Such a choice of numbers $\delta_p\ge 0$ is always possible since by log-convexity $p\mapsto\mu_p$ is non-decreasing and $(M_p)^{1/p}\rightarrow+\infty$ is equivalent to having $\lim_{p\rightarrow+\infty}\mu_p=+\infty\Leftrightarrow\sum_{j=1}^{+\infty}\delta_j=+\infty$ (e.g. see \cite[p.104]{compositionpaper}).

Conversely, given an arbitrary sequence $(\delta_p)_{p\ge 1}$ with $\delta_p\ge 0$ for all $p\ge 1$ and $\sum_{j=1}^{+\infty}\delta_j=+\infty$, then we can introduce a sequence $M\in\hyperlink{LCset}{\mathcal{LC}}$ via \eqref{muchoice}: We have $\mu_0=1$ (empty sum) and normalization follows by $\mu_1=\exp(\delta_1)\ge\exp(0)=1$.

The obtained sequence $M\in\hyperlink{LCset}{\mathcal{LC}}$ is unique since $\delta_1=\log(\mu_1)=\log(M_1/M_0)=\log(M_1)$ determines the first quotient. However, note that $(C^pM_p)_{p\in\NN}$, $C>0$ arbitrary, gives the quotients $(C\mu_p)_p$, hence yields the same sequence $(\delta_p)_{p\ge 2}$.\vspace{6pt}

In this notation property \hyperlink{mg}{$(\text{mg})$} holds if and only if
\begin{equation}\label{mgdelta}
\exists\;C\ge 1\;\forall\;p\in\NN:\;\;\;\log(\mu_{2p})-\log(\mu_p)=\sum_{l=p+1}^{2p}\delta_l\le C
\end{equation}
and \hyperlink{gamma1}{$(\gamma_1)$} resp. equivalently \hyperlink{beta1}{$(\beta_1)$} holds if and only if
\begin{equation}\label{snqdelta}
\exists\;Q\in\NN_{\ge 2}:\;\;\;\liminf_{p\rightarrow+\infty}\log(\mu_{Qp})-\log(\mu_p)=\liminf_{p\rightarrow+\infty}\sum_{l=p+1}^{Qp}\delta_l>\log(Q).
\end{equation}

Moreover we summarize:

\begin{remark}\label{mgremark}
\begin{itemize}
\item[$(i)$] \hyperlink{mg}{$(\on{mg})$} does imply both $\liminf_{p\rightarrow+\infty}\delta_p=0$ and $\sup_{p\ge 1}\delta_p<+\infty$.

If $\sup_{p\ge 1}\delta_p=+\infty$, then for all $k\in\NN$ we find some $p_k\in\NN$ with $\delta_{p_k}\ge k$ and so $\sum_{l=p_k}^{2p_k-2}\delta_l\ge k$ which makes \eqref{mgdelta} impossible. If $\liminf_{p\rightarrow+\infty}\delta_p>0$, then for some $\varepsilon>0$ (fixed) and $p_{\varepsilon}\in\NN$ we have $\delta_p\ge\varepsilon$ for all $p\ge p_{\varepsilon}$. This implies $\sum_{l=p+1}^{2p}\delta_l\ge\varepsilon p$ for all $p+1\ge p_{\varepsilon}$ making also \eqref{mgdelta} impossible as $p\rightarrow+\infty$.

\item[$(ii)$] However, the converse implication in $(i)$ is not true. Let $\delta_1:=1$ and $\delta_p:=\frac{1}{j+1}$ for $2^j+1\le p\le 2^{j+1}$, $j\in\NN$. Then $\delta_p\le 1$ for all $p\ge 1$ and $\delta_p\rightarrow 0$ as $p\rightarrow+\infty$, but $\sum_{l=2^j+1}^{2^{j+1}}\delta_l=\frac{2^j}{j+1}\rightarrow+\infty$ as $j\rightarrow+\infty$.
\end{itemize}
\end{remark}

Using this sequence $(\delta_p)_p$ we can prove new representations for the expressions $A_M(a_j,a_{j+1})$ and $B_M(a_j,a_{j+1})$.

\begin{lemma}\label{alternativedelta}
Let $(a_j)_{j\in\NN_{>0}}$ be an arbitrary sequence of positive integers with $a_{j+1}\ge a_j+2$ for all $j\ge 1$. Then we get
\begin{equation}\label{Aalternative}
A_M(a_j,a_{j+1})=\exp\left((a_{j+1}-a_j)\delta_{a_{j+1}+1}+\sum_{l=1}^{a_{j+1}-a_j-1}l\delta_{a_j+1+l}\right),
\end{equation}
and
\begin{equation}\label{Balternative}
B_M(a_j,a_{j+1})=\exp\left((a_{j+1}-a_j)\sum_{l=a_j+2}^{a_{j+1}}\delta_l-\sum_{l=1}^{a_{j+1}-a_j-1}l\delta_{a_j+1+l}\right).
\end{equation}
\end{lemma}

\demo{Proof}
We have
$$(\mu_{a_{j+1}+1})^{a_{j+1}-a_j}=\exp\left((a_{j+1}-a_j)\sum_{l=1}^{a_{j+1}+1}\delta_l\right),$$
and
\begin{align*}
\mu_{a_j+1}\cdots\mu_{a_{j+1}}&=\exp\left(\sum_{l=1}^{a_j+1}\delta_l+\sum_{l=1}^{a_j+2}\delta_l+\dots+\sum_{l=1}^{a_{j+1}}\delta_l\right)
\\&
=\exp\Bigg((a_{j+1}-a_j)\sum_{l=1}^{a_j+1}\delta_l+(a_{j+1}-a_j-1)\delta_{a_j+2}+(a_{j+1}-a_j-2)\delta_{a_j+3}
\\&
+\dots+(a_{j+1}-a_j-(a_{j+1}-a_j-1))\delta_{a_{j+1}}\Bigg)
\\&
=\exp\Bigg((a_{j+1}-a_j)\sum_{l=1}^{a_j+1}\delta_l+(a_{j+1}-a_j)\delta_{a_j+2}+\dots+(a_{j+1}-a_j)\delta_{a_{j+1}}
\\&
-\delta_{a_{j+2}}-2\delta_{a_{j+3}}-\dots-(a_{j+1}-a_j-1)\delta_{a_{j+1}}\Bigg)
\\&
=\exp\left((a_{j+1}-a_j)\sum_{l=1}^{a_j+1}\delta_l+(a_{j+1}-a_j)\sum_{l=1}^{a_{j+1}-a_j-1}\delta_{a_j+1+l}-\sum_{l=1}^{a_{j+1}-a_j-1}l\delta_{a_j+1+l}\right)
\\&
=\exp\left((a_{j+1}-a_j)\sum_{l=1}^{a_{j+1}}\delta_l-\sum_{l=1}^{a_{j+1}-a_j-1}l\delta_{a_j+1+l}\right).
\end{align*}

Similarly we have
$$(\mu_{a_j+1})^{a_{j+1}-a_j}=\exp\left((a_{j+1}-a_j)\sum_{l=1}^{a_j+1}\delta_l\right)$$
and so $$B_M(a_j,a_{j+1})=\exp\left((a_{j+1}-a_j)\sum_{l=1}^{a_{j+1}}\delta_l-\sum_{l=1}^{a_{j+1}-a_j-1}l\delta_{a_j+1+l}-(a_{j+1}-a_j)\sum_{l=1}^{a_j+1}\delta_l\right).$$
\qed\enddemo

By using \eqref{Aalternative} and \eqref{Balternative} we see in a better and more precise way that and how the size resp. growth of the sequences $(a_j)_j$, the ''Lusky numbers'', and $(\delta_j)_j$ is connected. Recall that by Lemma \ref{whatcanhappen} and Corollary \ref{whatcanhappencor} not each $M\in\hyperlink{LCset}{\mathcal{LC}}$ admits a sequence $(a_j)_j$ satisfying \eqref{22}, see also $(i)$ in Lemma \ref{alternativedeltaconsequ} below.

But, on the other hand when starting with a sequence $(a_j)_j$ satisfying some necessary growth restrictions, then we can show the following.

\begin{proposition}\label{ajexample}
Let $(a_j)_{j\ge 1}$ be a sequence of positive integers such that
\begin{equation}\label{Luskynecessary}
\forall\;j\ge 1:\;\;\;a_{j+1}-a_j\ge 2,\hspace{30pt}\sum_{j=1}^{+\infty}\frac{1}{a_{j+1}-a_j}=+\infty.
\end{equation}
Then there does exist a sequence $M\in\hyperlink{LCset}{\mathcal{LC}}$ such that the regularity condition $(b)$, i.e. \eqref{22}, holds true for this particular given $(a_j)_j$.
\end{proposition}

Recall that $a_{j+1}-a_j\ge 2$ is a necessary condition to have \eqref{22} for $(a_j)_j$ (see $(i)$ in Remark \ref{firstremark}). In Lemma \ref{alternativedeltaconsequ1} and Corollary \ref{alternativedeltaconsequ1cor} we show that also $\sum_{j=1}^{+\infty}\frac{1}{a_{j+1}-a_j}=+\infty$ is necessary for any $(a_j)_j$ to be considered in \eqref{22}. Consequently, this result tells us that for each sequence of integers $(a_j)_j$ which could be possibly used for Lusky numbers in the weight sequence setting there does exist a sequence $M\in\hyperlink{LCset}{\mathcal{LC}}$ such that $(a_j)_j$ are Lusky numbers for any weight $v_{M,c}$.

\demo{Proof}
We introduce $M$ by the sequence $(\delta_j)_j$ as follows:
$$\delta_{a_{j+1}+1}=\delta_{a_j+2}:=d_j,\;\;\;j\in\NN_{>0},\hspace{30pt}\delta_j:=0\;\;\;\text{else},$$
with
$$d_j:=\frac{C}{a_{j+1}-a_j+1},\hspace{30pt}3\le C<+\infty.$$
This growth shall be compared with $(ii)$ in Lemma \ref{alternativedeltaconsequ} below and note that $0<\frac{C}{a_{j+1}-a_j+1}\le\frac{C}{3}$. We get $\frac{1}{a_{j+1}-a_j+1}\ge\frac{1}{2}\frac{1}{a_{j+1}-a_j}\Leftrightarrow a_{j+1}-a_j\ge 1$, hence
$$\sum_{j=1}^{+\infty}\delta_j=\sum_{l=1}^{+\infty}\sum_{j=a_l+2}^{a_{l+1}+1}\delta_j=\sum_{l=1}^{+\infty}2d_l=2C\sum_{l=1}^{+\infty}\frac{1}{a_{l+1}-a_l+1}\ge C\sum_{l=1}^{+\infty}\frac{1}{a_{l+1}-a_l}=+\infty,$$
which shows $\lim_{p\rightarrow+\infty}\mu_p=+\infty$, so $M\in\hyperlink{LCset}{\mathcal{LC}}$ is verified.

By \eqref{Aalternative} and \eqref{Balternative} we get for all $j\ge 1$
$$A_M(a_j,a_{j+1})=\exp((a_{j+1}-a_j+1)d_j)=\exp(C)>2,\hspace{30pt}B_M(a_j,a_{j+1})=\exp((a_{j+1}-a_j-1)d_j).$$
Moreover $(a_{j+1}-a_j-1)d_j\le(a_{j+1}-a_j+1)d_j=C$ and $(a_{j+1}-a_j-1)d_j\ge 1\Leftrightarrow\frac{C}{a_{j+1}-a_j+1}\ge\frac{1}{a_{j+1}-a_j-1}\Leftrightarrow(C-1)(a_{j+1}-a_j)\ge C+1\Leftrightarrow a_{j+1}-a_j\ge\frac{C+1}{C-1}$ which holds true since $a_{j+1}-a_j\ge 2$ and $\frac{C+1}{C-1}\le 2\Leftrightarrow 3\le C$. Thus $B_M(a_j,a_{j+1})\ge\exp(1)>2$ holds true and so \eqref{22} is valid for $(a_j)_j$.
\qed\enddemo

Without any further information on the growth of $(a_j)_j$ it seems to be not possible to obtain further information on $M$; e.g. for having \hyperlink{mg}{$(\on{mg})$} we have to assume that $\sup_{j\ge 1}a_{j+1}-a_j<+\infty\Leftrightarrow\inf_{j\ge 1}\frac{1}{a_{j+1}-a_j}>0$ (see $(i)$ in \ref{mgremark} and also $(ii)$ in Lemma \ref{alternativedeltaconsequ}).

\begin{lemma}\label{alternativedeltaconsequ}
Let $(\delta_p)_{p\ge 1}$ be a sequence with $\delta_p\ge 0$ for all $p$ and $\sum_{j=1}^{+\infty}\delta_j=+\infty$. Let $M$ be the sequence obtained via \eqref{muchoice}.
\begin{itemize}
\item[$(i)$] If $\sup_{p\ge 1}\delta_p=+\infty$, which contradicts \hyperlink{mg}{$(\on{mg})$} for $M$, then there does not exist a strictly increasing sequence (of positive integers) $(a_j)_j$ satisfying \eqref{22} for $M$, see also the example in Lemma \ref{whatcanhappen}.

\item[$(ii)$] If $\limsup_{p\rightarrow+\infty}\delta_p<+\infty$, then each sequence $(a_j)_j$ enjoying \eqref{22} (for $M$) has to satisfy
    $$\sup_{j\ge 1}\delta_{a_{j+1}+1}(a_{j+1}-a_j)<+\infty,$$
    thus the growth rate of $(\delta_p)_p$ is limiting the maximal admissible growth of $(a_j)_j$.

In particular, if $\liminf_{p\rightarrow+\infty}\delta_p>0$, then $\sup_{j\ge 1}a_{j+1}-a_j<+\infty$, i.e.
$$\exists\;C\ge 2\;\forall\;j\ge 1:\;\;\;2\le a_{j+1}-a_j\le C,$$
or equivalently the space $H^{\infty}_{v_M,c}(\CC)$ is solid by \cite[Theorem 5.2]{BonetLuskyTaskinen19Bergman}, see Remark \ref{logsqureremark}.

\item[$(iii)$] If there exist numbers $0<d_1\le d_2<+\infty$ such that $d_1\le\delta_p\le d_2$ for all $p\in\NN_{>0}$, then \eqref{22} is satisfied for any sequence of integers $(a_j)_j$ satisfying $C_1\le a_{j+1}-a_j\le C_2$ for some $2\le C_1\le C_2$ and all $j\ge 1$ provided that
$$2\le d_1C_1(1+C_1),\hspace{30pt}1+\frac{d_2C_2(C_2-1)}{2}\le d_1C_1(C_1-1).$$
\end{itemize}
\end{lemma}

Roughly speaking this result shows that too fast increasing weight sequences do not admit the existence of Lusky numbers.

\demo{Proof}
$(i)$ By assumption we can find for any $k\in\NN$ some $p_k\in\NN$ such that $\delta_{p_k}\ge k$. Take now $(a_j)_{j\ge 1}$, then by $a_j\rightarrow+\infty$ we find $j_k\in\NN$ such that $a_{j_k}\le p_k\le a_{j_k+1}-1$. Because $a_{j+1}-a_j-1\ge 1$ has to be satisfied (see $(i)$ in Remark \ref{firstremark}) and since $\delta_p\ge 0$, by using \eqref{Aalternative} we can estimate for $A_M(a_j,a_{j+1})$ as follows:

If $p_k=a_{j_k}$, then $$A_M(a_{j_k-1},a_{j_k})\ge\exp\left(\sum_{l=1}^{a_{j_{k}}-a_{j_k-1}-1}l\delta_{a_{j_k-1}+1+l}\right)\ge\exp(\delta_{p_k}(a_{j_k}-a_{j_k-1}-1))\ge\exp(\delta_{p_k})\ge\exp(k).$$
If $a_{j_k}+2\le p_k\le a_{j_k+1}-1$, then similarly
$$A_M(a_{j_k},a_{j_k+1})\ge\exp\left(\sum_{l=1}^{a_{j_{k}+1}-a_{j_k}-1}l\delta_{a_{j_k}+1+l}\right)\ge\exp(\delta_{p_k}(p_k-a_{j_k}-1))\ge\exp(\delta_{p_k})\ge\exp(k).$$
Finally, if $p_k=a_{j_k}+1$, then
$$A_M(a_{j_k-1},a_{j_k})\ge\exp((a_{j_k}-a_{j_k-1})\delta_{p_k})\ge\exp(\delta_k)\ge\exp(k).$$

In any case, as $k\rightarrow+\infty$ we see that the upper bound in \eqref{22} becomes impossible. Note: As seen in $(i)$ in Remark \ref{mgremark}, the assumption $\sup_{p\ge 1}\delta_p=+\infty$ violates \hyperlink{mg}{$(\text{mg})$}.\vspace{6pt}

$(ii)$ This follows immediately by \eqref{Aalternative} because $A_M(a_j,a_{j+1})\ge\exp\left((a_{j+1}-a_j)\delta_{a_{j+1}+1}\right)$, hence $\sup_{j\ge 1}\delta_{a_{j+1}+1}(a_{j+1}-a_j)=+\infty$ violates the upper bound in \eqref{22}.

If $\liminf_{p\rightarrow+\infty}\delta_p>0$, then we can find some $\varepsilon>0$ small (from now on fixed) and $p_{\varepsilon}\in\NN$ such that $\delta_p\ge\varepsilon$ for all $p\ge p_{\varepsilon}$. In this situation, if $\sup_{j\ge 1}a_{j+1}-a_j=+\infty$ is valid, then for all $k\in\NN$ we can find $j_k\in\NN$, $j_k\ge p_{\varepsilon}$, such that $a_{j_k+1}-a_{j_k}\ge k$. Since $a_{j+1}-a_j\ge 1$ for all $j\ge 1$ we have $a_{j_k}\ge j_k\ge p_{\varepsilon}$ as well and so $\delta_{a_{j_k+1}+1}(a_{j_k+1}-a_{j_k})\ge\varepsilon k$ making $\sup_{j\ge 1}\delta_{a_{j+1}+1}(a_{j+1}-a_j)<+\infty$ impossible as $k\rightarrow+\infty$.\vspace{6pt}


$(iii)$ By assumption and \eqref{Aalternative} for each $j\ge 1$ we can estimate $A_M(a_j,a_{j+1})\le\exp(d_2C_2+d_2\frac{(a_{j+1}-a_j)(a_{j+1}-a_j-1)}{2})\le\exp(d_2C_2+d_2C_2^2)$, by \eqref{Balternative} we get $B_M(a_j,a_{j+1})\le\exp(d_2C_2(a_{j+1}-a_j-1))\le\exp(d_2C_2^2)$ (hence for the upper bounds we only require $d_2,C_2<+\infty$).

Moreover, for the estimate from below, it is enough to prove that the arguments in the exponents are bounded from below by $1$. We have that $A_M(a_j,a_{j+1})\ge\exp(C_1d_1+d_1\frac{(a_j-a_{j-1})(a_j-a_{j-1}-1)}{2})$, so it is required that $C_1d_1+\frac{d_1C_1(C_1-1)}{2}=\frac{2C_1d_1+d_1C_1^2-d_1C_1}{2}\ge 1\Leftrightarrow d_1C_1(C_1+1)\ge 2$. Moreover, for $B_M(a_j,a_{j+1})\ge\exp(1)$ we require that $(a_{j+1}-a_j)\sum_{l=a_j+2}^{a_{j+1}}\delta_l\ge 1+\sum_{l=1}^{a_{j+1}-a_j-1}l\delta_{a_j+1+l}$. Since $(a_{j+1}-a_j)\sum_{l=a_j+2}^{a_{j+1}}\delta_l\ge(a_{j+1}-a_j)(a_{j+1}-a_j-1)d_1\ge C_1(C_1-1)d_1$ and $\sum_{l=1}^{a_{j+1}-a_j-1}l\delta_{a_j+1+l}\le d_2\frac{(a_{j+1}-a_j-1)(a_{j+1}-a_j)}{2}\le d_2\frac{C_2(C_2-1)}{2}$ it is sufficient to require $d_1C_1(C_1-1)\ge 1+d_2\frac{C_2(C_2-1)}{2}$.
\qed\enddemo
	
\begin{example}\label{qGevreyexamples}
\begin{itemize}
\item[$(i)$] Consider $N^{q,\alpha}=(N^{q,\alpha}_p)_{p\in\NN}$ with $N^{q,\alpha}_p:=q^{p^{\alpha}}$, $q>1$ and $\alpha>2$. By \cite[Lemma 5.7]{compositionpaper} and \cite[Sect. 5.5]{whitneyextensionweightmatrix} we get that each $\omega_{N^{q,\alpha}}$ is equivalent (w.r.t. \hyperlink{sim}{$\sim$}) to the weight $\max\{0,\log(t)^s\}$ with $\alpha-1=\frac{1}{s-1}\Leftrightarrow s=\frac{1}{\alpha-1}+1<2$ (and $s>1$).

One has $\mu_p=q^{p^{\alpha}-(p-1)^{\alpha}}$ and so $\frac{\mu_p}{\mu_{p-1}}=q^{p^{\alpha}-2(p-1)^{\alpha}+(p-2)^{\alpha}}$, $p\ge 2$. Thus $\delta_p\rightarrow+\infty$ as $p\rightarrow+\infty$ follows and $(i)$ in Lemma \ref{alternativedeltaconsequ} implies that there does not exist a sequence $(a_j)_j$ satisfying \eqref{22}. This consequence should be compared with \cite[Corollary 2.8]{BonetTaskinen18}.

\item[$(ii)$] We apply $(iii)$ in Lemma \ref{alternativedeltaconsequ} to the $q$-Gevrey sequence $M^q=(M^q_p)_{p\in\NN}$ which yields $\delta_p=\log(\mu_p/\mu_{p-1})=\log(q^2)=2\log(q)=d_1=d_2$ for all $p\ge 1$. So $2\le d_1C_1(1+C_1)\Leftrightarrow 2\le 2\log(q)C_1(1+C_1)\Leftrightarrow e^{1/(C_1(C_1-1))}\le q$ has to be satisfied and $1+\frac{d_2C_2(C_2-1)}{2}\le d_1C_1(C_1-1)\Leftrightarrow 1\le 2\log(q)(C_1(C_1-1)-C_2(C_2-1)/2)$. Choosing $C_1=C_2$ gives $1\le\log(q)(C_1(C_1-1))$ and so again $e^{1/(C_1(C_1-1))}\le q$ is required. Consequently we can choose $C_1\in\NN_{\ge 2}$ such that $e^{1/(C_1(C_1-1))}\le q$ which should be compared with Proposition \ref{qGevreyprop} before.
\end{itemize}
\end{example}

Finally, we are going to prove an upper growth restriction for the sequence $(a_j)_j$ showing that also the second assumption in \eqref{Luskynecessary} above is really necessary.

\begin{lemma}\label{alternativedeltaconsequ1}
For any given $M\in\hyperlink{LCset}{\mathcal{LC}}$ a sequence (of positive integers) $(a_j)_{j\ge 1}$ enjoying \eqref{22} has to satisfy $\sum_{j=1}^{+\infty}\frac{1}{a_{j+1}-a_j}=+\infty$, which is ''restricting the growth of $(a_j)_j$ from above''.
\end{lemma}

The reason for writing ''restricting the growth of $(a_j)_j$ from above'' is that the divergence of the series is not excluding the situation having a subsequence $(j_k)_{k\ge 1}$ such that $a_{j_k+1}-a_{j_k}\ge c_k$ for all $k\ge 1$, with values $c_k\ge 1$ as large as desired.

\demo{Proof}
If $\sup_{p\ge 1}\delta_p=+\infty$, then there does not exist such a sequence $(a_j)_{j\ge 1}$, see $(i)$ in Lemma \ref{alternativedeltaconsequ}. Consequently, we have $\sup_{p\ge 1}\delta_p<+\infty$ and by $(ii)$ in Lemma \ref{alternativedeltaconsequ} it follows that $\delta_{a_{j+1}+1}\le\frac{C}{a_{j+1}-a_j}$ for some $C\ge 1$ and all $j\in\NN_{>0}$.

Assume now that $\sum_{j=1}^{+\infty}\frac{1}{a_{j+1}-a_j}<+\infty$ is valid for a sequence $(a_j)_j$ satisfying \eqref{22}. By \eqref{Aalternative} we get that $\sum_{l=1}^{a_{j+1}-a_j-1}l\delta_{a_j+1+l}\le C_1$ for some $C_1\ge 1$ and all $j\in\NN_{>0}$. This together with \eqref{Balternative} implies $$\sum_{l=a_j+2}^{a_{j+1}}\delta_l\le\frac{C_2}{a_{j+1}-a_j}+\frac{1}{a_{j+1}-a_j}\sum_{l=1}^{a_{j+1}-a_j-1}l\delta_{a_j+1+l}\le\frac{C_2}{a_{j+1}-a_j}+\frac{C_1}{a_{j+1}-a_j}$$ for some $C_2\ge 1$ and all $j\in\NN_{>0}$. Note that in each such sum we consider indices $\delta _i$ with $a_j+2\le i\le a_{j+1}$.

Hence we can estimate as follows:
\begin{align*}
\sum_{j=1}^{+\infty}\delta_j&=C_3+\sum_{j=a_1+2}^{+\infty}\delta_j=C_3+\sum_{j=1}^{+\infty}\delta_{a_{j+1}+1}+\sum_{j=1}^{+\infty}\sum_{l=a_j+2}^{a_{j+1}}\delta_{l}
\\&
\le C_3+C\sum_{j=1}^{+\infty}\frac{1}{a_{j+1}-a_j}+(C_1+C_2)\sum_{j=1}^{+\infty}\frac{1}{a_{j+1}-a_j}<+\infty.
\end{align*}
But this implies $\lim_{p\rightarrow+\infty}\mu_p<+\infty$, hence $\lim_{p\rightarrow+\infty}(M_p)^{1/p}<+\infty$, a contradiction to $M\in\hyperlink{LCset}{\mathcal{LC}}$.
\qed\enddemo

We summarize $(i)$ in Remark \ref{firstremark} and Lemma \ref{alternativedeltaconsequ1}:

\begin{corollary}\label{alternativedeltaconsequ1cor}
Let $M\in\hyperlink{LCset}{\mathcal{LC}}$ be given. Then any sequence (of positive integers) $(a_j)_{j\ge 1}$ enjoying \eqref{22} has to satisfy \eqref{Luskynecessary}.
\end{corollary}

Using this new information we can also prove now the following result which shows that even ''very regular and nice'' sequences considered in the ultradifferentiable and ultraholomorphic setting do not admit automatically a sequence $(a_j)_j$ satisfying \eqref{22}. (By ''nice'' we mean that the ultradifferentiable resp. ultraholomorphic function classes do satisfy several good stability properties.)

\begin{lemma}\label{whatcanhappen1}
There does exist $M\in\hyperlink{LCset}{\mathcal{LC}}$ having \hyperlink{mg}{$(\on{mg})$} and \hyperlink{beta1}{$(\beta_1)$} (i.e. $M$ is equivalent to $N\in\hyperlink{SRset}{\mathcal{SR}}$) but such that there does not exist a sequence (of positive integers) $(a_j)_{j\ge 1}$ satisfying \eqref{22}.
\end{lemma}

\demo{Proof}
Let $M$ be defined by its sequence of quotients $(\mu_p)_{p\ge 1}$ as follows: We put $\mu_0:=1$ and
$$\mu_p:=c_i,\;\;\;2^i\le p<2^{i+1},\;i\in\NN,$$
with $(c_i)_i$ a sequence of strictly increasing positive real numbers such that
\begin{itemize}
\item[$(i)$] $c_0\ge 1$,

\item[$(ii)$] $c_i\rightarrow+\infty$ as $i\rightarrow+\infty$,

\item[$(iii)$] $\liminf_{i\rightarrow+\infty}\frac{c_{i+1}}{c_i}>2$ and $\sup_{i\in\NN}\frac{c_{i+1}}{c_i}<+\infty$.

\end{itemize}
A straight-forward choice would be $c_i:=Q^i$, with $2<Q<+\infty$.\vspace{6pt}

{\itshape Claim:} $M\in\hyperlink{LCset}{\mathcal{LC}}$ and $M$ has both \hyperlink{mg}{$(\on{mg})$} and \hyperlink{beta1}{$(\beta_1)$}.

By definition $M\in\hyperlink{LCset}{\mathcal{LC}}$ is clear and moreover $\sup_{p\in\NN_{>0}}\frac{\mu_{2p}}{\mu_p}=\sup_{i\in\NN}\frac{c_{i+1}}{c_i}<+\infty$, which proves \hyperlink{mg}{$(\on{mg})$}, and $\liminf_{p\rightarrow+\infty}\frac{\mu_{2p}}{\mu_p}=\liminf_{i\rightarrow+\infty}\frac{c_{i+1}}{c_i}>2$, which proves \hyperlink{beta1}{$(\beta_1)$} (with $Q=2$ there).\vspace{6pt}

{\itshape Claim:} There does not exist a sequence (of integers) $(a_j)_{j\ge 1}$ satisfying \eqref{22}.

For convenience we set $I_l:=[2^l,2^{l+1})$, $l\in\NN$.

First we see that by definition and \eqref{Binteresting} we never can have $a_j,a_{j+1}\in I_l$ with $l\in\NN$ and $j\in\NN_{>0}$ arbitrary: In this situation $B_M(a_j,a_{j+1})=1$ follows, hence contradicting the lower estimate in \eqref{22}.

Let now $j\in\NN_{>0}$ be arbitrary (but fixed) and so $a_j\in I_{l_j}$ with some $l_j\in\NN_{>0}$. Then $a_{j+1}\in I_{l_j+i}$ with some $i\in\NN_{>0}$ follows and by assumption $(iii)$ we can assume that $\frac{c_{i+1}}{c_i}>2$ for all $i\ge l_j$. We have to distinguish between two cases:\vspace{6pt}

{\itshape Case I:} First, if $a_j+1=2^{l_j+1}$ and so $a_j+1\in I_{l_j+1}$, then $\mu_{a_j+1}=c_{l_j+1}$ and in order to avoid $B_M(a_j,a_{j+1})=1$ we have to choose $a_{j+1}\in I_{l_j+i}$, $i\ge 2$. Thus $a_{j+1}\ge 2^{l_j+2}$ follows which yields $a_{j+1}-a_j\ge 2^{l_j+2}-2^{l_j+1}+1=2^{l_j+1}+1$.

{\itshape Case II:} Second, if $a_j+1<2^{l_j+1}$, then $a_{j+1}\in I_{l_j+1}$ has to be valid for all $j$ large: If $a_{j+1}\in I_{l_j+i}$ with $i\ge 2$, then we would have by \eqref{Binteresting} that
$$B_M(a_j,a_{j+1})=\frac{\mu_{a_j+1}\mu_{a_j+2}\cdots\mu_{a_{j+1}}}{(\mu_{a_j+1})^{a_{j+1}-a_j}}\ge\left(\frac{c_{l_j+1}}{c_{l_j}}\right)^{2^{l_j+2}-2^{l_j+1}}>2^{2^{l_j+1}},$$
which tends to infinity as $j\rightarrow+\infty$. Thus the upper estimate in \eqref{22} fails for large $j$. Similarly we see that in Case I above we have $a_{j+1}\in I_{l_j+2}$ because $i\ge 3$ would imply $B_M(a_j,a_{j+1})\ge\left(\frac{c_{l_j+2}}{c_{l_j+1}}\right)^{2^{l_j+3}-2^{l_j+2}}>2^{2^{l_j+2}}$, again contradicting the upper estimate in \eqref{22} as $j\rightarrow+\infty$.\vspace{6pt}

On the other hand, in the first case when $a_j+1=2^{l_j+1}$ and $a_{j+1}\in I_{l_j+2}$, then again by \eqref{Binteresting} we have
$$B_M(a_j,a_{j+1})\ge\left(\frac{c_{l_j+2}}{c_{l_j+1}}\right)^{a_{j+1}-2^{l_j+2}+1}>2^{a_{j+1}-2^{l_j+2}+1},$$
and in the second one, when $a_j+1<2^{l_j+1}$ and $a_{j+1}\in I_{l_j+1}$, then we get
$$B_M(a_j,a_{j+1})\ge\left(\frac{c_{l_j+1}}{c_{l_j}}\right)^{a_{j+1}-2^{l_j+1}+1}>2^{a_{j+1}-2^{l_j+1}+1}.$$
Consequently, in order to guarantee the upper estimate in \eqref{22} necessary requirements are $a_{j+1}-2^{l_j+2}\le d$ and $a_{j+1}-2^{l_j+1}\le d$ for some $d\in\NN_{>0}$ not depending on $j\ge 1$.

Finally we recall that in the first case, as shown above, the choice $a_{j+1},a_{j+2}\in I_{l_j+2}$ is not possible and so this restriction implies $a_{j+2}-a_{j+1}\ge 2^{l_j+3}-2^{l_j+2}-d=2^{l_j+2}-d$. In the second case, similarly the choice $a_{j+1},a_{j+2}\in I_{l_j+1}$ is not possible which implies $a_{j+2}-a_{j+1}\ge 2^{l_j+2}-2^{l_j+1}-d=2^{l_j+1}-d$.\vspace{6pt}

We summarize: Eventually only the second case can occur and so for all $i\in\NN_{>0}$ large enough we have $a_{j+i}\in I_{l_j+i}$, $a_{j+i+1}\in I_{l_j+i+1}$ and $a_{j+i}\le 2^{j+i}+d$, hence $a_{j+i+1}-a_{j+i}\ge 2^{l_j+i}-d\rightarrow+\infty$ as $i\rightarrow+\infty$. But this would imply $\sum_{k\ge 1}\frac{1}{a_{k+1}-a_k}<+\infty$ which contradicts Corollary \ref{alternativedeltaconsequ1cor}.



\qed\enddemo

\vspace{6pt}
We close this section with the following observations for $M\in\hyperlink{LCset}{\mathcal{LC}}$.

\begin{itemize}
\item[$(i)$] If there does exist a sequence (of integers) $(a_j)_{j\ge 1}$ such that
$$\exists\;1<b\le 2\;\forall\;j\in\NN_{>0}:\;\;\;b\le\min\{A_M(a_j,a_{j+1}),B_M(a_j,a_{j+1})\},$$
then by iteration we can get the first estimate in \eqref{22}: One has
\begin{align*}
A_M(a_j,a_{j+2})&:=\frac{(\mu_{a_{j+2}+1})^{a_{j+2}-a_j}}{\mu_{a_{j+2}}\cdots\mu_{a_j+1}}=\frac{(\mu_{a_{j+2}+1})^{a_{j+2}-a_{j+1}}}{\mu_{a_{j+2}}\cdots\mu_{a_{j+1}+1}}\cdot\frac{(\mu_{a_{j+2}+1})^{a_{j+1}-a_j}}{\mu_{a_{j+1}}\cdots\mu_{a_j+1}}
\\&
\ge\frac{(\mu_{a_{j+2}+1})^{a_{j+2}-a_{j+1}}}{\mu_{a_{j+2}}\cdots\mu_{a_{j+1}+1}}\cdot\frac{(\mu_{a_{j+1}+1})^{a_{j+1}-a_j}}{\mu_{a_{j+1}}\cdots\mu_{a_j+1}}=A_M(a_{j+1},a_{j+2})\cdot A_M(a_j,a_{j+1})\ge b^2,
\end{align*}
and similarly
$$B_M(a_j,a_{j+2}):=\frac{\mu_{a_j+1}\cdots\mu_{a_{j+2}}}{(\mu_{a_j+1})^{a_{j+2}-a_j}}\ge\frac{\mu_{a_j+1}\cdots\mu_{a_{j+1}}}{(\mu_{a_j+1})^{a_{j+1}-a_j}}\cdot\frac{\mu_{a_{j+1}+1}\cdots\mu_{a_{j+2}}}{(\mu_{a_{j+1}+1})^{a_{j+2}-a_{j+1}}}=B_M(a_j,a_{j+1})\cdot B_M(a_{j+1},a_{j+2}).$$
Iterating these estimates $n$ times, $n\in\NN_{>0}$ chosen minimal such that $b^n>2$, we are done. Consequently, the choice $a'_j:=a_{nj}$ yields a sequence satisfying the first estimate in \eqref{22}. However, if
$$\exists\;K\ge 1\;\forall\;j\in\NN_{>0}:\;\;\max\{A_M(a_j,a_{j+1}),B_M(a_j,a_{j+1})\}\le K,$$
then in general it is not clear that the upper estimate in \eqref{22} holds true for the sequence $a'_j:=a_{nj}$.

\item[$(ii)$] Such an iteration, yielding a stretching of the Lusky numbers, is always possible if for \eqref{22} any sequence $(a_j)_j$ satisfying $C_1\le a_{j+1}-a_j\le C_2$ for some $2\le C_1\le C_2$ can be used (e.g. for the $q$-Gevrey sequences, see Proposition \ref{qGevreyprop}): In this case $a'_j:=a_{nj}$ satisfies $nC_1\le a'_{j+1}-a'_j\le nC_2$ for all $j\ge 1$ and is still ''admissible''.

\item[$(iii)$] But in general, if $(a_j)_{j\ge 1}$ is a strictly increasing sequence of positive integers satisfying \eqref{22}, then an arbitrary subsequence $a'_k:=(a_{j_k})_{k\ge 1}$ cannot be considered for satisfying \eqref{22}. More precisely, given such $(a_j)_{j\ge 1}$, then $a'_j:=a_{2^j}$ can never be considered: By using the necessary assumption $2\le a_{j+1}-a_j$ we obtain
    $$a'_{j+1}-a'_j=a_{2^{j+1}}-a_{2^j}=\left(a_{2^{j+1}}-a_{2^{j+1}-1}\right)+\dots+\left(a_{2^j+1}-a_{2^j}\right)\ge 2^j\cdot 2=2^{j+1},$$
which yields $\sum_{j\ge 1}\frac{1}{a'_{j+1}-a'_j}<+\infty$, hence a contradiction by Corollary \ref{alternativedeltaconsequ1cor}. (Note that in $(ii)$ above $\sum_{j\ge 1}\frac{1}{a'_{j+1}-a'_j}=+\infty$ is clear.)
\end{itemize}

\subsection{Second example - Gevrey sequences}\label{Gevreyexample}
We use the technique from the previous section and consider for $s>0$ and $j\in\NN_{>0}$ the sequence $\delta_j:=\frac{s}{j}$.

\begin{proposition}\label{harmonicseriesprop}
Let $M^s\in\hyperlink{LCset}{\mathcal{LC}}$ be defined by its sequence of quotients $\mu^s_p$ via $\delta_j:=\frac{s}{j}$ and \eqref{muchoice}, i.e. $\mu^s_p:=\exp\left(s\sum_{l=1}^{p}\frac{1}{l}\right)$, $p\ge 1$, and $\mu^s_0:=1$.

By the asymptotic growth behavior of $p\mapsto\sum_{l=1}^{p}\frac{1}{l}$ it follows that $C^{-1}p^s\le\mu^s_p\le Cp^s$ for some $C\ge 1$ and all $p\in\NN$, hence $M^s$ is equivalent to the Gevrey sequence $G^s=(p!^s)_{p\in\NN}$.\vspace{6pt}

Then \eqref{22} is satisfied for the sequence (of integers) $a_j:=c(j+5)^2$, $c:=\lceil\frac{1}{s}\rceil$.

In particular, if $s\ge 1$, then \eqref{22} is satisfied for the sequence (of integers) $a_j:=(j+5)^2$.
\end{proposition}

The equivalence between $M^s$ and $G^s$ does imply (see Corollary \ref{hullsequiv})
$$\forall\;c,d>0:\;\;\;S(H^{\infty}_{v_{M^s,c}}(\CC))\cong S(H^{\infty}_{v_{G^s,d}}(\CC)),\hspace{30pt}s(H^{\infty}_{v_{M^s,c}}(\CC))\cong s(H^{\infty}_{v_{G^s,d}}(\CC)),$$
and $\omega_{M^s}(t)=O(\omega_{G^s}(Ct))$ and $\omega_{G^s}(t)=O(\omega_{M^s}(Ct))$ for some $C\ge 1$ as $t\rightarrow+\infty$. Since each power weight $t\mapsto at^{1/s}$, $s,a>0$ arbitrary, does have $\omega(2t)=O(\omega(t))$ (which is a standard assumption in the theory of ultradifferentialbe functions, see \cite{BraunMeiseTaylor90}) we get that $\omega_{M^s}\hyperlink{sim}{\sim}\omega_{G^s}\hyperlink{sim}{\sim}t\mapsto at^{1/s}$ for each $a>0$. Moreover, each arising weight does have \eqref{omega6} as well.

Altogether we have obtained: $\forall\;a>0\;\forall\;b,c,d>0$:
$$S(H^{\infty}_{v_{M^s,b}}(\CC))\cong S(H^{\infty}_{v_{G^s,c}}(\CC))\cong S(H^{\infty}_{w^{1/s}_{d,a}}(\CC)),\;\;\;s(H^{\infty}_{v_{M^s,b}}(\CC))\cong s(H^{\infty}_{v_{G^s,c}}(\CC))\cong s(H^{\infty}_{w^{1/s}_{d,a}}(\CC)),$$
with $w^{1/s}_{d,a}$ denoting the weight $t\mapsto\exp(-a(dt)^{1/s})$. Thus Proposition \ref{harmonicseriesprop} gives an alternative proof for the representation obtained in \cite[Theorem 3.1]{BonetTaskinen18} (see also Remark \ref{forwardshift}).

\demo{Proof}
Obviously we have $cs\ge 1$ and moreover $cs\le(\frac{1}{s}+1)s=1+s$ does hold true.

Then \eqref{Aalternative} turns into $$A_M(a_j,a_{j+1})=\exp\left((a_{j+1}-a_j)\frac{s}{a_{j+1}+1}+s\sum_{l=1}^{a_{j+1}-a_j-1}\frac{l}{a_j+1+l}\right)$$
and \eqref{Balternative} turns into $$B_M(a_j,a_{j+1})=\exp\left((a_{j+1}-a_j)s\sum_{l=a_j+2}^{a_{j+1}}\frac{1}{l}-s\sum_{l=1}^{a_{j+1}-a_j-1}\frac{l}{a_j+1+l}\right).$$

In order to simplify the computations and since the first three claims below hold true for $a_j:=c(j+2)^2$ we use now this sequence in the following computations. First, we have $a_{j+1}-a_j=c(j+3)^2-c(j+2)^2=c(2j+5)$, $j\ge 1$, and start with $A_M(a_j,a_{j+1})$.\vspace{6pt}

{\itshape Claim I:} $A_M(a_j,a_{j+1})\ge\exp(1)>2$ for all $j\ge 1$. For any sequence of integers $(a_j)_j$ we have $(a_{j+1}-a_j)\frac{s}{a_{j+1}+1}\ge 0$ (which tends to $0$ as $j\rightarrow+\infty$ by our choice $a_j$) and
$$s\sum_{l=1}^{a_{j+1}-a_j-1}\frac{l}{a_j+1+l}\ge s\frac{1}{a_{j+1}}\sum_{l=1}^{a_{j+1}-a_j-1}l=s\frac{(a_{j+1}-a_j-1)(a_{j+1}-a_j)}{2a_{j+1}}.$$
Hence we require
\begin{align*}
&s\sum_{l=1}^{a_{j+1}-a_j-1}\frac{l}{a_j+1+l}\ge s\frac{(c(2j+5)-1)c(2j+5)}{2c(j+3)^2}\ge 1,
\end{align*}
i.e. $j^2(4cs-2)+j(20cs-2s-12)+25sc-5s-18\ge 0$ which holds true for all $j\ge 1$: We have $cs\ge 1$ and so $25cs=5cs+20cs\ge 5s+20$. Moreover, for $0<s\le 1$ we get $20cs-2s-12\ge 20-12=8$ and for $s>1$ we get $c=1$, hence $20cs-2s-12=18s-12>6$.\vspace{6pt}

{\itshape Claim II:} $A_M(a_j,a_{j+1})\le\exp(8cs)\le\exp(8(1+s))$ for all $j\ge 1$. The first arising summand tends to $0$ as $j\rightarrow+\infty$, more precisely we have $(a_{j+1}-a_j)\frac{s}{a_{j+1}+1}=cs(2j+5)\frac{1}{c(j+3)^2+1}\le A\Leftrightarrow 2csj+5cs\le Acj^2+A6cj+A9c+A$ for all $j\ge 1$ by choosing $A:=\frac{5cs}{9c+1}$ because $A6c\ge 2cs\Leftrightarrow\frac{15c}{9c+1}\ge 1\Leftrightarrow 6c\ge 1$. Note that $A\le 4cs\Leftrightarrow 5\le 36c+4\Leftrightarrow 1\le 36c$.

For the second summand we estimate as follows:
\begin{align*}
&s\sum_{l=1}^{a_{j+1}-a_j-1}\frac{l}{a_j+1+l}\le s\frac{1}{a_j+2}\sum_{l=1}^{a_{j+1}-a_j-1}l=s\frac{(a_{j+1}-a_j-1)(a_{j+1}-a_j)}{2(a_j+2)}=s\frac{(c(2j+5)-1)c(2j+5)}{2(c(j+2)^2+2)}
\\&
=s\frac{4c^2j^2+20c^2j-2cj+25c^2-5c}{2cj^2+8cj+8c+4}\le A\Longleftrightarrow 4c^2j^2+jc(20c-2)+25c^2-5c\le\frac{A}{s}(2cj^2+8cj+8c+4),
\end{align*}
with $A:=4sc$, which yields the claim.\vspace{6pt}

{\itshape Claim III:} $B_M(a_j,a_{j+1})\le\exp(8cs)\le\exp(8(1+s))$ for all $j\ge 1$. Since the second summand is always $\le 0$ we have to study the first one. We estimate as follows:
\begin{align*}
&(a_{j+1}-a_j)s\sum_{l=a_j+2}^{a_{j+1}}\frac{1}{l}\le(a_{j+1}-a_j)s\frac{a_{j+1}-a_j-1}{a_j+2}\le A_1,
\end{align*}
for all $j\ge 1$ with the choice $A_1:=8sc=2A$ by Claim II above.\vspace{6pt}

{\itshape Claim IV:} $B_M(a_j,a_{j+1})\ge\exp(1)>2$ for all $j\ge 4$, which explains the index shift.

We get
$$(a_{j+1}-a_j)s\sum_{l=a_j+2}^{a_{j+1}}\frac{1}{l}\ge s(a_{j+1}-a_j)(a_{j+1}-a_j-1)\frac{1}{a_{j+1}}$$
and
$$s\sum_{l=1}^{a_{j+1}-a_j-1}\frac{l}{a_j+1+l}\le s\frac{(a_{j+1}-a_j-1)(a_{j+1}-a_j)}{2(a_j+2)},$$
hence
\begin{align*}
B_M(a_j,a_{j+1})&\ge\exp\left(s(a_{j+1}-a_j)(a_{j+1}-a_j-1)\frac{1}{a_{j+1}}-s\frac{(a_{j+1}-a_j-1)(a_{j+1}-a_j)}{2(a_j+2)}\right)
\\&
=\exp\left(s(a_{j+1}-a_j-1)(a_{j+1}-a_j)\left(\frac{1}{a_{j+1}}-\frac{1}{2(a_j+2)}\right)\right)
\\&
=\exp\left(s(a_{j+1}-a_j-1)(a_{j+1}-a_j)\left(\frac{2(a_j+2)-a_{j+1}}{2a_{j+1}(a_j+2)}\right)\right).
\end{align*}
We show that the last argument arising in $\exp(\cdot)$ is $\ge 1$ for all $j\ge 1$. First,  $s(a_{j+1}-a_j-1)(a_{j+1}-a_j)=s(4c^2j^2+20c^2j-2cj+25c^2-5c)$ (see Claim II) and $2(a_j+2)-a_{j+1}=cj^2+2cj-c+4\ge 0$ for all $j\ge 1$. Moreover we have $s(4c^2j^2+20c^2j-2cj+25c^2-5c)(cj^2+2cj-c+4)\ge(4cj^2+20cj-2j+25c-5)(cj^2+2cj-c+4)$ and the difference between $(4cj^2+20cj-2j+25c-5)(cj^2+2cj-c+4)$ and $2a_{j+1}(a_j+2)=2c(j+3)^2(c(j+2)^2+2)$ is given by
$$2c^2j^4+8c^2j^3-2cj^3-13c^2j^2+3cj^2-90c^2j+48cj-8j-97c^2+69c-20\ge 0,$$
which holds true for all $j\ge 4$.
\qed\enddemo

\section{From weight functions to weight sequences}\label{fromfcttosequ}
The aim of this section is to see how the weight sequence setting is becoming meaningful when starting with an abstractly given weight function $v:[0,+\infty)\rightarrow(0,+\infty)$ in the weighted holomorphic setting, i.e. $v$ is continuous, non-increasing and rapidly decreasing. We call $v$ {\itshape normalized}, when $v(t)=1$ for all $t\in[0,1]$ which can be assumed w.l.o.g.: Otherwise replace $v$ by $w$ such that $w$ is normalized and $v(t)=w(t)$ for all $t\ge t_0>1$, which yields $H^{\infty}_v(\CC)=H^{\infty}_w(\CC)$.

\begin{lemma}\label{vBMTweight}
Let $v$ be a normalized weight function. Then
\begin{equation}\label{omegafromv}
\omega^v(t):=-\log(v(t)),\;\;\;t\in[0,+\infty),
\end{equation}
satisfies the following conditions arising frequently in the theory of ultradifferentiable classes defined by so-called {\itshape Braun-Meise-Taylor weight functions} (see \cite{BraunMeiseTaylor90}):
\begin{itemize}
\item[$(i)$] $\omega^v:[0,+\infty)\rightarrow[0,+\infty)$ is continuous, non-decreasing, $\lim_{t\rightarrow+\infty}\omega^v(t)=+\infty$ and $\omega^v(t)=0$ for $t\in[0,1]$ (normalization), i.e. $(\omega_0)$ in \cite[Sect. 2.2]{dissertation},
\item[$(ii)$] $\log(t)=o(\omega^v(t))$ as $t\rightarrow+\infty$, i.e. $(\omega_3)$ in \cite[Sect. 2.2]{dissertation}.
\end{itemize}
$v$ does satisfy in addition
\begin{equation}\label{vconvexity}
t\mapsto-\log(v(e^t))\;\;\;\text{is convex on}\;\RR,
\end{equation}
if and only if $\varphi_{\omega^v}: t\mapsto\omega^v(e^t)$ is convex, i.e. $\omega^v$ has $(\omega_4)$ in \cite[Sect. 2.2]{dissertation}.
\end{lemma}

For concrete given $v$, condition \eqref{vconvexity} may be checked by straight-forward computations, e.g. it holds true for the weight $v(t):=\exp(-\exp(t))$ mentioned in \cite[Rem. 2.2]{BonetLuskyTaskinen19} and \cite[Cor. 2.8]{BonetTaskinen18} and also for the weight $v(t)=\exp(-\max\{0,\log(t)^2\})$.

If $v\equiv v_M$ for $M\in\hyperlink{LCset}{\mathcal{LC}}$, and so $v$ is a normalized weight function because $\omega_M(t)=0$ for $t\in[0,1]$, then $\omega^v\equiv\omega_M$ and so $\varphi_{\omega^v}=\varphi_{\omega_M}$ is always convex by definition (e.g. see \cite[Lemma 12, $(4)\Rightarrow(5)$]{BonetMeiseMelikhov07}).

\demo{Proof}
$(i)$ follows immediately by definition and the properties for $v$.

$(ii)$ We have to show that for each $\varepsilon>0$ there does exist $D_{\varepsilon}\ge 1$ such that $\log(t)\le\varepsilon\omega^v(t)+D_{\varepsilon}$ for all $t>0$, which is equivalent to $t\le(v(t))^{-\varepsilon}\exp(D_{\varepsilon})$ and so to $v(t)t^{1/\varepsilon}\le\exp(D_{\varepsilon}/\varepsilon)$. This holds true because $v$ is rapidly decreasing.
\qed\enddemo

The next result shows that for given weight functions $v$ satisfying all requirements from before we can associate a weight sequence such that the corresponding weight functions $v$ and $v_M$ can be compared.

\begin{proposition}\label{vBMTweight1}
Let $v$ be a normalized weight function satisfying \eqref{vconvexity} and $\omega^v$ be the weight given by \eqref{omegafromv}. Then the associated weight sequence $M^v=(M^v_p)_{p\in\NN}$ defined by
\begin{equation}\label{vBMTweight1equ1}
M^v_p:=\sup_{t>0}\frac{t^p}{\exp(\omega^v(t))}=\sup_{t>0}t^pv(t),
\end{equation}
belongs to the set \hyperlink{LCset}{$\mathcal{LC}$} and we get $\omega_{M^v}\hyperlink{sim}{\sim}\omega^v$, more precisely $\exists\;A\ge 1\;\forall\;t\in[0,+\infty):$
\begin{equation}\label{vBMTweight1equ2}
\frac{1}{A}(v_{M^{v}}(t))^{2}=\frac{1}{A}\exp(-2\omega_{M^v}(t))\le v(t)=\exp(-\omega^v(t))\le\exp(-\omega_{M^v}(t))=v_{M^v}(t).
\end{equation}
\end{proposition}

\demo{Proof}
We recall the definition of the {\itshape Legendre-Fenchel-Young-conjugate} of $\varphi_{\omega^v}$ by
$$\varphi^{*}_{\omega^v}(x):=\sup\{x y-\varphi_{\omega^v}(y): y\ge 0\},\;\;\;x\ge 0.$$
Since $\varphi_{\omega^v}$ is non-decreasing, convex by assumption, $\varphi_{\omega^v}(0)=0$ (by normalization) and $\lim_{t\rightarrow+\infty}\frac{\varphi_{\omega^v}(t)}{t}=\lim_{s\rightarrow+\infty}\frac{\omega^v(s)}{\log(s)}=+\infty$ by $(\omega_3)$, we get: $\varphi^{*}_{\omega^v}$ is convex, $\varphi^{*}_{\omega^v}(0)=0$, $x\mapsto\frac{\varphi^{*}_{\omega^v}(x)}{x} $ is non-decreasing and tending to $+\infty$ as $x\rightarrow+\infty$, e.g. see \cite[Rem. 1.3]{BraunMeiseTaylor90}.

Hence for $p\in\NN$
\begin{align*}
M^v_p&=\sup_{t>0}\frac{t^p}{\exp(\omega^v(t))}=\exp\sup_{t>0}\left(p\log(t)-\omega^v(t)\right)=\exp\sup_{s\in\RR}\left(ps-\omega^v(e^s)\right)=\exp\sup_{s\ge0}\left(ps-\omega^v(e^s)\right)
\\&
=\exp(\varphi^{*}_{\omega^v}(p)),
\end{align*}
which implies $M^v\in\hyperlink{LCset}{\mathcal{LC}}$ by the properties of the conjugate and note that by normalization $\omega^v(e^s)=0$ for all $-\infty<s\le 0$.

The equivalence $\omega_{M^v}\hyperlink{sim}{\sim}\omega^v$ has been shown in \cite[Theorem 4.0.3]{dissertation}, see also \cite[Lemma 5.7]{compositionpaper}. Note that for this proof also property $(\omega_4)$ on $\omega^v$ has to be used and more precisely we have shown that
$$\exists\;C\ge 1\;\forall\;t\ge 0:\;\;\;\omega_{M^v}(t)\le\omega^v(t)\le 2\omega_{M^v}(t)+C,$$
which yields \eqref{vBMTweight1equ2}.
\qed\enddemo

In the case $v\equiv v_M$, i.e. $\omega^v\equiv\omega_M$, we get
$$M_p=\sup_{t>0}\frac{t^p}{\exp(\omega_M(t))}=\sup_{t>0}\frac{t^p}{\exp(\omega^v(t))}=M^v_p$$
for all $p\in\NN$, hence $M^v\equiv M$. The first equality holds by \cite[p. 17]{mandelbrojtbook} (see also \cite[Proposition 3.2]{Komatsu73}), which has motivated the definition \eqref{vBMTweight1equ1}.\vspace{6pt}

We close this section with the following observations:

By \eqref{vBMTweight1equ1}, from given $v$ we compute $M^v$ via $M^v_p=t_p^pv(t_p)$ with $t_p$ denoting the global maximum point of $t\mapsto t^pv(t)$ (recalling the notation from \cite[Sect. 2]{BonetTaskinen18} and \cite[Def. 2.1]{BonetLuskyTaskinen19}). For $p=0$ we have $M^v_0=1$ by normalization. We put $\mu^v_0:=1$ and moreover get
\begin{equation}\label{rmucomparison}
\forall\;p\in\NN:\;\;\;t_p\le\mu^v_{p+1}\le t_{p+1}.
\end{equation}
This holds because $\mu^v_{p+1}=\frac{M^v_{p+1}}{M^v_p}=\frac{(t_{p+1})^{p+1}v(t_{p+1})}{(t_p)^{p}v(t_{p})}$ and by definition, $(t_{p+1})^{p+1}v(t_{p+1})\ge t^{p+1}v(t)$ for all $t\ge 0$, in particular $(t_{p+1})^{p+1}v(t_{p+1})\ge t_p^{p+1}v(t_p)$ which proves $\mu^v_{p+1}\ge t_p$. Similarly $t_{p+1}t_p^pv(t_p)\ge t_{p+1}t^pv(t)$ for all $t\ge 0$ is valid, in particular $t_{p+1}t_p^pv(t_p)\ge (t_{p+1})^{p+1}v(t_{p+1})$ which proves the second half.\vspace{6pt}

For the expressions under consideration for the regularity condition $(b)$ for $v$ (see \cite[$(2.1)$]{BonetTaskinen18}), for any integers $0<k<l$ we get
$$A_v(k,l):=\left(\frac{t_k}{t_l}\right)^k\frac{v(t_k)}{v(t_l)}=\frac{t_k^kv(t_k)t_l^l}{t_l^lv(t_l)t_l^k}=\frac{M^v_k}{M^v_l}t_l^{l-k}=\frac{t_l^{l-k}}{\mu^v_{k+1}\cdots\mu^v_l},$$
$$B_v(k,l):=\left(\frac{t_l}{t_k}\right)^l\frac{v(t_l)}{v(t_k)}=\frac{t_l^lv(t_l)t_k^k}{t_k^kv(t_k)t_k^l}=\frac{M^v_l}{M^v_k}t_k^{k-l}=\frac{\mu^v_{k+1}\cdots\mu^v_l}{t_k^{l-k}},$$
which should be compared, by taking into account \eqref{rmucomparison}, with \eqref{Ainteresting} and \eqref{Binteresting} for $k=a_j$ and $l=a_{j+1}$ and in particular with Remark \ref{ABexpressionsremark} for $M\equiv M^v$.

Consequently, the study of and search for (integer) Lusky numbers for given normalized weight $v$ satisfying \eqref{vconvexity} does precisely correspond to the study of this problem for the associated weight sequence $M^v$.


\section{Solid hulls and cores of weighted holomorphic functions on disks}\label{disksection}
In \cite{BonetTaskinendisc18}, see the main Theorem 2.2 there, also the solid hull of weighted holomorphic functions defined on the unit disk has been computed for some exponential weights, which are corresponding to the Gevrey weights with index $s\ge\frac{1}{2}$ in our setting \eqref{radialweights} below. It has turned out that for the computations the regularity condition $(b)$ has become crucial as well by using also in this framework the proof of \cite[Theorem 2.5]{BonetTaskinen18}, see \cite[Theorem 2.1]{BonetTaskinendisc18}.

Finally, in \cite{BonetLuskyTaskinen19} the solid hull and core for more general weight functions has been characterized, again by using condition $(b)$ and including all Gevrey weights when taking the generalizing function $w\equiv 1$ in the main result \cite[Theorem 3.1]{BonetLuskyTaskinen19}, see also \cite[Example 3.3 $(i)$]{BonetLuskyTaskinen19}.

Unfortunately, in general in this setting when considering holomorphic functions on a (the unit) disk it seems to be much more technical and difficult to determine a sequence (of positive real numbers) satisfying condition $(b)$ and so to obtain concrete representations for the solid hulls and cores by involving the ''Lusky numbers''.

We see now that also for the weight sequence setting the arising expressions are becoming much more complicated than before.\vspace{6pt}

For any $M\in\hyperlink{LCset}{\mathcal{LC}}$ and $c>0$ we set $\mathbb{D}^c:=\{z\in\CC: |z|<\frac{1}{c}\}$ and
\begin{equation}\label{radialweights}
v_{M,\mathbb{D}^c}(r):=\exp\left(-\omega_M\left(\frac{1}{1-cr}\right)\right),\;\;\;r\in\left[0,\frac{1}{c}\right).
\end{equation}
We put $\mathbb{D}:=\mathbb{D}^1$, i.e. denoting the unit disk. Then we introduce spaces of holomorphic functions on the disc $\mathbb{D}^c$ as follows:
\begin{align*}
H^{\infty}_{v_{M,\mathbb{D}^c}}(\mathbb{D}^c)&:=\{f\in H(\mathbb{D}^c): \|f\|_{v_{M,\mathbb{D}^c}}:=\sup_{z\in\mathbb{D}^c}|f(z)|v_{M,\mathbb{D}^c}(|z|)<+\infty\}
\\&
=\{f\in H(\mathbb{D}^c): \|f\|_{v_{M,\mathbb{D}^c}}:=\sup_{z\in\mathbb{D}^c}|f(z)|\exp\left(-\omega_M\left(\frac{1}{1-c|z|}\right)\right)<+\infty\},
\end{align*}
with $H(\mathbb{D}^c)$ denoting the class of all holomorphic functions on $\mathbb{D}^c$.

For each $c>0$ the mapping $r\mapsto v_{M,\mathbb{D}^c}(r)$ is continuous, non-increasing and $\lim_{r\rightarrow 1/c}r^kv_{M,\mathbb{D}^c}(r)=0$ for all $k\ge 0$, i.e. rapidly decreasing and so a weight function on $\mathbb{D}^c$. By normalization we have $\mu_1\ge 1$, hence $\omega_M(1)=0$ and so $v_{M,\mathbb{D}^c}(0)=1$.

\begin{remark}\label{alldiskequiv}
One shall observe that $H^{\infty}_{v_{M,\mathbb{D}^c}}(\mathbb{D}^c)\cong H^{\infty}_{v_{M,\mathbb{D}}}(\mathbb{D})$ for any $c>0$ isometrically:

Given $c>0$ and an arbitrary weight $v$ defined on $\mathbb{D}^c$, then consider $w:\mathbb{D}\rightarrow(0,+\infty)$ given by $w(z):=v(z/c)$ for all $z\in\mathbb{D}$. The operator $T: H_w^{\infty}(\mathbb{D})\rightarrow H_v^{\infty}(\mathbb{D}^c)$ given by $T(f)(z):=f(zc)$, $z\in\mathbb{D}^c$, $f\in H_w^{\infty}(\mathbb{D})$, realizes the isometric isomorphism since $v(z)T(f)(z)=w(zc)f(zc)$, $z\in\mathbb{D}^c$, $f\in H_w^{\infty}(\mathbb{D})$.

The analogous proof yields that $H^{\infty}_{v_{M,c}}(\CC)\cong H^{\infty}_{v_{M}}(\CC)$ for any $c>0$ isometrically and this corresponds to Corollary \ref{hullsequiv} and the fact that the arising expressions for condition $(b)$ are not depending on given $c>0$.

Moreover, when $f(z)=\sum_{j=0}^{+\infty}a_jz^j$ and $T(f)(z)=\sum_{j=0}^{+\infty}b_jz^j$, then $a_j=\frac{b_j}{c^j}$ for all $j\in\NN$ follows and this shall be compared with $c$ arising in the denominator in the characterizations in Theorems \ref{Thm25} and \ref{Thm24}.\vspace{6pt}

However, also in this setting, the equivalence $\omega_M\hyperlink{sim}{\sim}\omega_N$ will in general not imply the equivalence of the classes $H^{\infty}_{v_{M,\mathbb{D}^c}}(\mathbb{D}^c)$ and also not imply that the corresponding solid hulls and cores are isomorphic.
\end{remark}

Apart from the power-weights corresponding to the Gevrey-sequences we are pointing out the following two examples:

First recall that, as commented in Sect. \ref{qGevrey} above for each $q$-Gevrey sequence $M^q=(q^{p^2})_{p\in\NN}$, $q>1$, we have that $\omega_{M^q}$ is equivalent to the weight $\max\{0,\log(t)^2\}$. Hence this case corresponds on the unit disk (up to equivalence of weight functions) to
$$v(r):=\exp(-(\log(\frac{1}{1-r}))^2)=\exp(-(\log(s))^2)=\frac{1}{s^{\log(s)}}=(1-r)^{-\log(1-r)},\,\;\;\;0\le r<1,$$
by taking $s:=\frac{1}{1-r}$ (and so $s\in[1,+\infty)$).\vspace{6pt}

Second, for given $\alpha>0$ on $\mathbb{D}$ we can consider the ''standard weight'' $v_{\alpha}(z):=(1-|z|)^{\alpha}$, see \cite[Remark 2.2]{BonetLuskyTaskinen19}. In this case we get
$$(1-r)^{\alpha}=\exp\left(-w\left(\frac{1}{1-r}\right)\right)\Longleftrightarrow-\alpha\log(1-r)=w\left(\frac{1}{1-r}\right)\Longleftrightarrow-\alpha\log\left(\frac{1}{s}\right)=w(s),\;\;\;s\in[1,+\infty),$$
hence $w(s)=\alpha\log(s)$ for some weight $w$ (taking again $s:=\frac{1}{1-r}$). Thus, in order to apply the weight sequence approach, by \eqref{radialweights} we are searching for some $M\in\hyperlink{LCset}{\mathcal{LC}}$ such that $\omega_M\hyperlink{sim}{\sim}t\mapsto\alpha\log(1+t)$, equivalently that $\omega_M\hyperlink{sim}{\sim}t\mapsto\log(1+t)$, $t\ge 0$. However, we show now that we cannot find such a sequence $M$, hence these ''standard weights'' on the unit disk cannot be considered within the set \hyperlink{LCset}{$\mathcal{LC}$}.\vspace{6pt}

\begin{lemma}\label{lognotallowed}
	Concerning the function $t\mapsto\log(1+t)$ on $[0,+\infty)$ we get:
\begin{itemize}
\item[$(i)$] There does not exist $M\in\hyperlink{LCset}{\mathcal{LC}}$ satisfying
\begin{equation}\label{logequi}
\omega_M\hyperlink{sim}{\sim}t\mapsto\log(1+t).
\end{equation}
\item[$(ii)$] Let $M=(M_p)_p$ be a sequence with $1=M_0$ satisfying
$$\exists\;q_0\in\NN_{>0}\;\forall\;p>q_0:\;\;\;M_p=+\infty,$$
and such that $1\le\mu_p\le\mu_{p+1}$ for $1\le p\le q_0$ with $\mu_{q_0+1}=\frac{M_{q_0+1}}{M_{q_0}}=+\infty$, i.e. $M_p\in\RR_{>0}$ and $M$ is log-convex for only finitely many indices $p$. Then $\omega_M$ does have \eqref{logequi}, when $\omega_M$ is defined via \eqref{assofunc} by using the conventions $\frac{1}{+\infty}=0$ and $\log(0)=-\infty$.
\end{itemize}
\end{lemma}

The first part means that any weight $a\log(1+t)$, $a>0$, cannot occur in the equivalence class of associated weight functions coming from (standard) weight sequences. However, the second part yields that for ''exotic sequences'', i.e. $M_p=+\infty$ for all but finitely many $p\in\NN$, each weight $a\log(1+t)$, $a>0$, is equivalent to $\omega_M$.

\demo{Proof}
$(i)$ Assume that there exists some $M\in\hyperlink{LCset}{\mathcal{LC}}$ satisfying \eqref{logequi}. Then $\omega_M(t)\le C\log(1+t)+C$ for some $C\ge 1$ and all $t\ge 0$ and by using \cite[Proposition 3.2]{Komatsu73} we get for all $p\in\NN$:
\begin{align*}
M_p=\sup_{t>0}\frac{t^p}{\exp(\omega_M(t))}\ge\frac{1}{e^C}\sup_{t>0}\frac{t^p}{\exp(C\log(1+t))}=\frac{1}{e^C}\left(\sup_{t>0}\frac{t^{p/C}}{1+t}\right)^C.
\end{align*}
But for all $p>C$ we see that $\frac{t^{p/C}}{1+t}\rightarrow+\infty$ as $t\rightarrow+\infty$ and so $M_p=+\infty$ for all $p>C$, a contradiction to having $M\in\RR_{>0}^{\NN}$.

$(ii)$ By definition and the conventions we have $\omega_M(0):=0$ and $\omega_M(t)=\sup_{p\in\NN_0}\log\left(\frac{t^p}{M_p}\right)=\max_{0\le p\le q_0}\log\left(\frac{t^p}{M_p}\right)$, $t>0$. Moreover, $\frac{t^p}{M_p}\le\frac{t^{p+1}}{M_{p+1}}\Leftrightarrow\mu_{p+1}\le t$, $p\ge 0$, is valid and so:

If $0<t\le\mu_1$, then $\omega_M(t)=0$ and if $\mu_{p_t}\le t<\mu_{p_t+1}$ for some $1\le p_t\le q_0-1$, then $\omega_M(t)=p_t\log(t)-\log(M_{p_t})$. Finally, if $\mu_{q_0}\le t<\mu_{q_0+1}=+\infty$ (i.e. for all $t$ large enough), then $\omega_M(t)=q_0\log(t)-\log(M_{q_0})$ follows which proves \eqref{logequi}.
\qed\enddemo

Now we proceed studying the weights $v_{M,\mathbb{D}^c}$ following the computations given in Sect. \ref{regularitybcond}. Let from now on $c>0$ be arbitrary but fixed.\vspace{6pt}

For any $k\ge 0$ and $r\ge 0$ we set $G^k_{M,\mathbb{D},c}(r):=r^kv_{M,\mathbb{D}^c}(r)$, hence we get via using the auxiliary function $h_M$ (recall $h_M(t):=\inf_{p\in\NN}t^pM_p$, $t>0$, $h_M(0):=0$, hence $h_M(t)=\exp(-\omega_M(1/t))$) that
\begin{align*}
G^k_{M,\mathbb{D},c}(r)&=r^k\exp\left(-\omega_M\left(\frac{1}{1-cr}\right)\right)=r^kh_M(1-cr),\;\;\;0\le r<\frac{1}{c}.
\end{align*}

Since $M\in\hyperlink{LCset}{\mathcal{LC}}$, by recalling \eqref{assovsmu} we have $h_M(1-cr)=1$ for all $0\le r\le\frac{1}{c}-\frac{1}{c\mu_1}$ (note that $\frac{1}{c}-\frac{1}{c\mu_1}\ge 0\Leftrightarrow\mu_1\ge 1$ by normalization) and for $p\ge 1$:
$$h_M(1-cr)=(1-cr)^pM_p\;\;\;\text{if}\;\frac{1}{\mu_{p+1}}\le 1-cr<\frac{1}{\mu_{p}}\Longleftrightarrow\frac{1}{c}-\frac{1}{c\mu_p}<r\le\frac{1}{c}-\frac{1}{c\mu_{p+1}}.$$
Since $h_M$ is continuous, also $G^k_{M,\mathbb{D},c}$ is so and for convenience we set now $I_{p,c}:=(\frac{1}{c}-\frac{1}{c\mu_p},\frac{1}{c}-\frac{1}{c\mu_{p+1}}]$, $p\ge 1$, and $I_{0,c}:=[0,\frac{1}{c}-\frac{1}{c\mu_1}]$.

With $r_{k,\mathbb{D},c}$ we denote the global maximum point of $G^k_{M,\mathbb{D},c}$.

If $k=0$, then $G^k_{M,\mathbb{D},c}(r)=h_M(1-cr)$ and so each $r\in I_{0,c}$ is a maximum point (there we have $h_M\equiv 1$ and $h_M$ is non-decreasing).

Let now $k>0$ be fixed, then on $I_{0,c}$ we have $G^k_{M,\mathbb{D},c}(r)=r^k$, which is strictly increasing with maximum (right end) point $\frac{1}{c}-\frac{1}{c\mu_1}$.

For all $p\in\NN$ and $r\in[0,\frac{1}{c})$ we set $f_{k,p,c}(r):=r^k(1-cr)^pM_p$, i.e. $f_{k,p,c}\equiv G^k_{M,\mathbb{D},c}|_{I_{p,c}}$. We see that $f_{k,p,c}(r)\le f_{k,p+1,c}(r)\Leftrightarrow r\le\frac{1}{c}-\frac{1}{c\mu_{p+1}}$, in particular this is valid for all $r\in I_{p,c}$.

A direct computation for $p\in\NN_{>0}$ yields $f'_{k,p,c}(r)=kr^{k-1}(1-cr)^pM_p+r^kp(1-cr)^{p-1}(-c)M_p$, hence the unique maximum point of $f_{k,p,c}$ on $[0,\frac{1}{c})$ is given by $s_{k,p,c}:=\frac{k}{c(k+p)}$ (for $p=0$ we would have for $r\mapsto r^k$ that such a maximum point does not exist anymore because $\frac{1}{c}$ does not belong to the interval $[0,\frac{1}{c})$).

Note that $0<s_{k,p,c}<\frac{1}{c}$ is valid because $k<k+p$. Moreover we have $\frac{k}{c(k+1)}=s_{k,1,c}>s_{k,p,c}>s_{k,p+1,c}$, $\lim_{k\rightarrow+\infty}s_{k,p,c}=\frac{1}{c}$ for all $p\in\NN_{>0}$ and $\lim_{p\rightarrow+\infty}s_{k,p,c}=0$ for all $k>0$.\vspace{6pt}

We have to distinguish now and summarize:

\begin{itemize}
\item[$(i)$] For all small $k>0$ satisfying $s_{k,1,c}=\frac{k}{c(k+1)}\le\frac{1}{c}-\frac{1}{c\mu_1}\Leftrightarrow k\le\mu_1-1$, hence $s_{k,p,c}\in I_{0,c}$ for all $p\ge 1$ as well holds true, each function $f_{k,p,c}$ is (strictly) decreasing on $I_{p,c}$, $p\ge 1$. Thus $r_{k,\mathbb{D},c}=\frac{1}{c}-\frac{1}{c\mu_1}$ for all such small $k>0$. (Recall that by normalization we have $\mu_1\ge 1$.)

\item[$(ii)$] For all larger values $k>0$ we have that $s_{k,q,c}\in I_{p,c}$ for at least one pair $p,q\ge 1$.

If $s_{k,q,c}\in I_{p,c}$ with $q>p\ge 1$, then $s_{k,q,c}<s_{k,p,c}$ and $f_{k,p,c}$ is (strictly) increasing at $s_{k,q,c}$. Since $f_{k,p,c}\equiv G^k_{M,\mathbb{D},c}|_{I_{p,c}}$ we get that also $G^k_{M,\mathbb{D},c}$ is strictly increasing at $s_{k,q,c}$. Similarly, if $1\le q<p$, then $f_{k,p,c}$ is (strictly) decreasing at $s_{k,q,c}$ and so $G^k_{M,\mathbb{D},c}$ as well.

Finally, if $s_{k,p,c}\in I_{p,c}$ for some $p\ge 1$, then $s_{k,p,c}$ is the maximum point for $f_{k,p,c}$ on $[0,\frac{1}{c})$ resp. for $G^k_{M,\mathbb{D},c}$ on $I_{p,c}$.

\item[$(iii)$] Note that we always find at least one such value $p\ge 1$ because $\lim_{p\rightarrow+\infty}s_{k,p,c}=0$ and $\lim_{p\rightarrow+\infty}\frac{1}{c}-\frac{1}{c\mu_p}=\frac{1}{c}$. Moreover, when found such integer $p$, then we have $s_{k,q,c}<s_{k,p,c}<s_{k,t,c}$ for all $q>p>t\ge 1$. Consequently, by the definition of the intervals $I_{p,c}$ and because $f_{k,p,c}\equiv G^k_{M,\mathbb{D},c}|_{I_{p,c}}$ we see that $r_{k,\mathbb{D},c}=s_{k,p,c}$ is valid.
\end{itemize}

Thus we have shown for the maximum value point $r_{k,\mathbb{D},c}$ of $G^k_{M,\mathbb{D},c}$, $k\ge 0$:
\begin{equation}\label{maximumvaluefordiscs}
r_{k,\mathbb{D},c}=\frac{1}{c}-\frac{1}{c\mu_1},\;\;\;\text{for}\;0\le k\le\mu_1-1,\hspace{15pt}r_{k,\mathbb{D},c}=s_{k,p_k,c}=\frac{k}{c(k+p_k)},\;\text{s. th.}\;s_{k,p_k,c}\in I_{p_k,c},\;\;\;\text{for}\;k>\mu_1-1.
\end{equation}
Thus, for any given $k>\mu_1-1$, we have $\frac{1}{c}-\frac{1}{c\mu_{p_k}}<r_{k,\mathbb{D},c}\le\frac{1}{c}-\frac{1}{c\mu_{p_k+1}}$ for some $p_k\ge 1$.

Now assume w.l.o.g. that $\mu_2>1$ and then, for any $p\in\NN_{>0}$, we have
\begin{align*}
\frac{k}{c(k+p)}&=\frac{1}{c}-\frac{1}{c\mu_{p_k+1}}\Longleftrightarrow\frac{k}{(k+p)}=\frac{\mu_{p+1}-1}{\mu_{p+1}}\Longleftrightarrow\frac{\mu_{p+1}}{\mu_{p+1}-1}=1+\frac{k}{p}
\\&
\Longleftrightarrow\frac{1}{\mu_{p+1}-1}=\frac{p}{k}\Longleftrightarrow k=p(\mu_{p+1}-1).
\end{align*}
Thus, when choosing $k_p:=p(\mu_{p+1}-1)$ and so $k_p<k_{p+1}$ for all $p\ge 1$, $k_p\rightarrow+\infty$ as $p\rightarrow+\infty$, we obtain
$$r_{k_p,\mathbb{D},c}=\frac{k_p}{c(k_p+p)}=\frac{1}{c}-\frac{1}{c\mu_{p+1}},$$
which implies by \cite[1.8. III]{mandelbrojtbook} that $$v_{M,\mathbb{D}^c}(r_{k_p,\mathbb{D},c})=\exp\left(-\omega_M\left(\frac{1}{1-cr_{k_p,\mathbb{D},c}}\right)\right)=\exp(-\omega_M(\mu_{p+1}))=\frac{M_{p+1}}{(\mu_{p+1})^{p+1}}.$$
Next, for arbitrary $1\le p<q$ (with $p,q\in\NN$) we introduce the expressions
$$A_{M,\mathbb{D}^c}(p,q):=\left(\frac{r_{k_p,\mathbb{D},c}}{r_{k_q,\mathbb{D},c}}\right)^{k_p}\frac{v_{M,\mathbb{D}^c}(r_{k_p,\mathbb{D},c})}{v_{M,\mathbb{D}^c}(r_{k_q,\mathbb{D},c})},\hspace{30pt}B_{M,\mathbb{D}^c}(p,q):=\left(\frac{r_{k_q,\mathbb{D},c}}{r_{k_p,\mathbb{D},c}}\right)^{k_q}\frac{v_{M,\mathbb{D}^c}(r_{k_q,\mathbb{D},c})}{v_{M,\mathbb{D}^c}(r_{k_p,\mathbb{D},c})},$$
and so we can prove:

\begin{lemma}\label{ABexpressionsdisc}
Let $M\in\hyperlink{LCset}{\mathcal{LC}}$ be given (with having $1\le\mu_1<\mu_2$) and $p,q\in\NN_{>0}$ with $1\le p<q$. Then for any $c>0$ we have
$$A_{M,\mathbb{D}^c}(p,q)=\left(\frac{\mu_{p+1}-1}{\mu_{q+1}-1}\right)^{p(\mu_{p+1}-1)}\cdot\frac{(\mu_{q+1})^{p(\mu_{p+1}-1)+q}}{(\mu_{p+1})^{p\mu_{p+1}+1}\mu_{p+2}\cdots\mu_q},$$
$$B_{M,\mathbb{D}^c}(p,q)=\left(\frac{\mu_{q+1}-1}{\mu_{p+1}-1}\right)^{q(\mu_{q+1}-1)}\cdot\frac{(\mu_{p+1})^{q(\mu_{q+1}-1)+p}\mu_{p+2}\cdots\mu_q}{(\mu_{q+1})^{q\mu_{q+1}}}.$$
\end{lemma}

\demo{Proof}
By the computations before we get $\left(\frac{r_{k_p,\mathbb{D},c}}{r_{k_q,\mathbb{D},c}}\right)^{k_p}=\left(\frac{1-\frac{1}{\mu_{p+1}}}{1-\frac{1}{\mu_{q+1}}}\right)^{p(\mu_{p+1}-1)}=\left(\frac{(\mu_{p+1}-1)\mu_{q+1}}{(\mu_{q+1}-1)\mu_{p+1}}\right)^{p(\mu_{p+1}-1)}$ and $\frac{v_{M,\mathbb{D}^c}(r_{k_p,\mathbb{D},c})}{v_{M,\mathbb{D}^c}(r_{k_q,\mathbb{D},c})}=\frac{M_{p+1}}{M_{q+1}}\frac{(\mu_{q+1})^{q+1}}{(\mu_{p+1})^{p+1}}=\frac{(\mu_{q+1})^{q+1}}{\mu_{p+2}\cdots\mu_{q+1}(\mu_{p+1})^{p+1}}=\frac{(\mu_{q+1})^q}{\mu_{p+1}\mu_{p+2}\cdots\mu_q(\mu_{p+1})^p}$. Putting this information together we are done, the expression for $B_{M,\mathbb{D}^c}(p,q)$ follows analogously.
\qed\enddemo

Like in the weighted entire function case this result shows that the expressions are not depending on the parameter $c>0$ anymore, hence we can write $A_{M,\mathbb{D}}(p,q)$ and $B_{M,\mathbb{D}}(p,q)$ instead and the letter $\mathbb{D}$ shall emphasize the different behavior compared with the entire case in Section \ref{solidhullsandcores}. Finally, applying \cite[Theorems 2.3, 2.4]{BonetLuskyTaskinen19} and Remark \eqref{alldiskequiv}, we obtain the following result:

\begin{theorem}\label{Thm2524disc}
Let $M\in\hyperlink{LCset}{\mathcal{LC}}$ be given (with having $1\le\mu_1<\mu_2$) such that there exists a strictly increasing sequence (of integers) $(a_j)_{j\in\NN_{\ge 1}}$ (''Lusky numbers'') and numbers $b$ and $K$ with $K\ge b>2$ such that
\begin{equation}\label{Thm2524discequ}
b\le\min\{A_{M,\mathbb{D}}(a_j,a_{j+1}),B_{M,\mathbb{D}}(a_j,a_{j+1})\}\le\max\{A_{M,\mathbb{D}}(a_j,a_{j+1}),B_{M,\mathbb{D}}(a_j,a_{j+1})\}\le K,
\end{equation}
i.e. satisfying the regularity condition $(b)$. Then the solid hull of $H^{\infty}_{v_{M,\mathbb{D}^c}}(\mathbb{D}^c)$ is given by
\begin{equation*}\label{Thm2524discequ1}
	S(H^{\infty}_{v_{M,\mathbb{D}^c}}(\mathbb{D}^c))=\left\{(b_j)_{j\in\NN}\in\CC^{\NN}: \sup_{j\in\NN_{\ge 1}}v_{M,\mathbb{D}^c}(r_{k_{a_j},\mathbb{D},c})\left(\sum_{a_j<l\le a_{j+1}}|b_l|^2(r_{k_{a_j},\mathbb{D},c})^{2l}\right)^{1/2}<+\infty\right\},
\end{equation*}
or equivalently by
\begin{equation*}\label{Thm2524discequ2}
 S(H^{\infty}_{v_{M,\mathbb{D}^c}}(\mathbb{D}^c))=\left\{(b_j)_{j\in\NN}\in\CC^{\NN}: \sup_{j\in\NN_{\ge 1}}\frac{M_{a_j+1}}{(\mu_{a_j+1})^{a_j+1}}\left(\sum_{a_j<l\le a_{j+1}}|b_l|^2\left(\frac{1}{c}-\frac{1}{c\mu_{a_j+1}}\right)^{2l}\right)^{1/2}<+\infty\right\}.
\end{equation*}
The solid core of $H^{\infty}_{v_{M,\mathbb{D}^c}}(\mathbb{D}^c)$ is given by
\begin{equation*}\label{Thm2524discequ3}
	s(H^{\infty}_{v_{M,\mathbb{D}^c}}(\mathbb{D}^c))=\left\{(b_j)_{j\in\NN}\in\CC^{\NN}: \sup_{j\in\NN_{\ge 1}}v_{M,\mathbb{D}^c}(r_{k_{a_j},\mathbb{D},c})\sum_{a_j<l\le a_{j+1}}|b_l|(r_{k_{a_j},\mathbb{D},c})^{l}<+\infty\right\},
\end{equation*}
	or equivalently by
\begin{equation*}\label{Thm2524discequ4}
	s(H^{\infty}_{v_{M,\mathbb{D}^c}}(\mathbb{D}^c))=\left\{(b_j)_{j\in\NN}\in\CC^{\NN}: \sup_{j\in\NN_{\ge 1}}\frac{M_{a_j+1}}{(\mu_{a_j+1})^{a_j+1}}\sum_{a_j<l\le a_{j+1}}|b_l|\left(\frac{1}{c}-\frac{1}{c\mu_{a_j+1}}\right)^{l}<+\infty\right\}.
\end{equation*}
\end{theorem}

By these representations it is immediate to see again like in the entire case above that all solid hulls and cores are isomorphic w.r.t. the arising parameter $c>0$, i.e. not depending on the radius of the disk. The problem is now again to find a sequence $(a_j)_{j\ge 1}$ satisfying \eqref{Thm2524discequ}, which is in this setting much more complicated since one has to treat the expressions obtained in Lemma \ref{ABexpressionsdisc}.

\bibliographystyle{plain}
\bibliography{Bibliography}

\end{document}